\documentclass[a4paper]{amsart}

\usepackage[T1]{fontenc}
\usepackage{lmodern}

\usepackage{mathtools,amsthm,amssymb}

\newcommand{\E}   {\operatorname{E}}
\newcommand{\va}  {\operatorname{var}}
\newcommand{\bias}{\operatorname{bias}}
\newcommand{\err} {\operatorname{rmse}}
\newcommand{\cost}{\operatorname{cost}}
\newcommand{\tV}  {{\widetilde{V}}}
\newcommand{\tX}  {{\widetilde{X}}}
\newcommand{\tY}  {{\widetilde{Y}}}
\newcommand{\eps} {\varepsilon}
\newcommand{\bit} {\mathrm{bit}}
\newcommand{\num} {\mathrm{c}}
\newcommand{\hh}  {\mathrm{tMil}}

\newcommand{\R}{{\mathbb{R}}}
\newcommand{\N}{{\mathbb{N}}}

\DeclareMathOperator{\var}{Var}

\usepackage{tikz}
\usepackage{caption}
\usepackage{subcaption}
\usepackage{graphicx}
\usepackage{pgfplots}
\pgfplotsset{compat=1.14}
\usepackage{pgfgantt}
\usepackage{pdflscape}

\theoremstyle{plain}
\newtheorem{lemma}{Lemma}

\begin{document}

\title[A Random Bit Multilevel Algorithm]
{An Adaptive Random Bit Multilevel Algorithm for SDEs}

\author[Giles]
{Michael B.\ Giles}
\address{Mathematical Institute\\
University of Oxford\\
Oxford OX2 6GG\\
England} 
\email{mike.giles@maths.ox.ac.uk} 

\author[Hefter]
{Mario Hefter}
\address{Fachbereich Mathematik\\
Technische Universit\"at Kaisers\-lautern\\
Postfach 3049\\
67653 Kaiserslautern\\
Germany}
\email{\{hefter,lmayer,ritter\}@mathematik.uni-kl.de}

\author[Mayer]
{Lukas Mayer}

\author[Ritter]
{Klaus Ritter}

\begin{abstract}
We study the approximation of expectations $\E(f(X))$
for solutions $X$ of stochastic differential equations
and functionals $f$ on the path space
by means of Monte Carlo algorithms that only use random bits instead 
of random numbers. We construct an adaptive random 
bit multilevel algorithm, which is based on the Euler scheme,
the L\'evy-Ciesielski representation 
of the Brownian motion,
and asymptotically optimal random bit
approximations of the standard normal distribution.
We numerically compare this algorithm with
the adaptive classical multilevel Euler algorithm for a
geometric Brownian motion, an Ornstein-Uhlenbeck process,
and a Cox-Ingersoll-Ross process.
\end{abstract}


\date{February 28, 2019}

\maketitle

\section{Introduction}

We study the approximation of expectations $\E(f(X))$, where 
$X=(X(t))_{t\in [0,1]}$ is the $r$-dimensional solution of an autonomous 
stochastic differential equation (SDE) driven by a $d$-dimensional 
Brownian motion and where $f \colon C([0,1],\R^r)  \to \R$
is a functional on the path space. The main contribution 
of this paper is the construction of an adaptive random bit multilevel 
algorithm $A^\bit_\eps$, which is based on the generic adaptive multilevel 
algorithm from \cite{MLMC}. Here $\eps>0$ is an accuracy demand and 
input to the algorithm, and the maximal level as well as the replication 
numbers per level are determined adaptively. 
For a survey on multilevel Monte Carlo algorithms we refer to \cite{G15}. 

The algorithm $A^\bit_\eps$ employs 
the Euler scheme, the L\'evy-Ciesielski representation of the 
Brownian motion (Brownian bridge construction), 
and the asymptotically optimal random bit approximation 
of the standard normal distribution according to \cite[Thm.~1]{GHMR17}. 
Unfortunately, we have no analysis 
for the error and the cost of $A^\bit_\eps$, and even a proof of the 
convergence $\lim_{\eps \to 0} A^\bit_\eps(f) = \E(f(X))$ (under
suitable assumptions on the coefficients of the SDE and the
functional $f$) is missing. Instead we present numerical experiments. 

On each level $\ell$ a multilevel algorithm has to couple a fine
approximation and a course approximation. In the classical setting,
where random numbers are used, one may simply simulate Brownian increments 
for $2^\ell + 1$ equidistant points. For the fine approximation the Euler 
scheme with $2^\ell$ steps is applied, and for the course approximation 
the step-size is doubled and the corresponding increments are added up.
Of course, there are several options for the simulation of the
Brownian increments, in particular, one may either
simulate the increments directly or
use the L\'evy-Ciesielski representation of the Brownian motion.

Both of these approaches may be adapted to the random bit setting
by approximating all the involved normally distributed random variables 
by random variables that can be simulated with random bits only. However, 
in contrast to the classical setting, the two constructions no longer 
end up in the same distribution. The first approach, where increments 
are approximated and the approximations are added up, has been presented 
and analyzed in \cite{GHMR19}. 
Here the random bit approximations of the increments are independent, 
but an additional bias is introduced, in contrast to the classical
setting. 

In this paper we present the second approach,
which has already been sketched in \cite[Sec.~10.2]{G15}, and
which employs the L\'evy-Ciesielski representation and random
bit approximations to its normally distributed coefficients.
Here the approximations of the Brownian increments are no longer 
independent, which forms an obstacle for an error analysis,
but we obtain matching distributions:
The distribution of the course approximation on level $\ell$
coincides with the distribution of the fine approximation on level
$\ell-1$. Furthermore, this approach is well suited as the
building block for an adaptive multilevel algorithm.

In the numerical experiments we apply the adaptive random bit
multilevel Euler algorithm $A^\bit_\eps$ and the 
adaptive classical multilevel Euler algorithm $A^\num_\eps$,
which is based on random numbers, for three processes
$X$ and functionals $f$, namely, the maximum of a geometric Brownian
motion and the terminal values of an Ornstein-Uhlenbeck process and 
of a Cox-Ingersoll-Ross process. 
At first we compare the building blocks, i.e., the 
random bit Euler scheme and the classical Euler scheme, in terms
of their bias and variance decays. The decays depend on the process
$X$ and the functional $f$ under consideration, but in all three cases 
we observe no essential difference between the random bit and the 
classical Euler scheme. Next we turn to the adaptive algorithms.
In all three cases and for both algorithms the actual root mean squared 
error is almost equal to the accuracy demand $\eps$.
Finally, to achieve the same root mean squared error
the number of random bits needed by $A^\bit_\eps$ is only about $4$ to $7$ 
times larger than the number of random numbers needed by 
$A^\num_\eps$.

For the terminal value of the Cox-Ingersoll-Ross process we also
apply the truncated Milstein scheme from \cite{MR3732573}
as the building block, instead of the Euler scheme.
This leads to a substantially faster decay of the variance in the
classical and in the random bit case. Moreover,
the algorithm based on the truncated Milstein scheme and random bits
performs as good as the algorithm based on the Euler scheme and random 
numbers.

\section{Euler Schemes}\label{s2}

Consider an autonomous system 
\[
\phantom{\qquad\quad t \in [0,1]}
\mathrm{d}X(t) = a(X(t))\,\mathrm{d}t + b(X(t))\,\mathrm{d}W(t), 
\qquad\quad t \in [0,1],
\]
of SDEs with a deterministic initial value $X(0) = x_0 \in \R^r$ and a 
$d$-dimensional Brownian motion $W$, and with drift and diffusion 
coefficients $a \colon \R^r \to \R^r$ and 
$b \colon \R^r \to \R^{r \times d}$, respectively. Furthermore, consider 
the time discretization given by 
\begin{equation}\label{g1}
\phantom{\qquad\quad k=0,\dots,2^\ell,}
t_{k,\ell} = k /2^\ell,
\qquad\quad k=0,\dots,2^\ell,
\end{equation}
together with a suitable choice of
\begin{equation}\label{g5}
V_\ell = (V_{1,\ell}, \dots, V_{2^\ell,\ell}) 
\end{equation}
with $d$-dimensional random vectors $V_{k,\ell}$ on a common
probability space. These random vectors are meant to at least
approximate the Brownian increments associated to \eqref{g1}, and
the corresponding Euler scheme $X_\ell$ is given by
$X_\ell(t_{0,\ell}) = x_0$ and
\begin{align}\label{euler1}
X_\ell(t_{k,\ell}) = X_\ell(t_{k-1,\ell}) + 2^{-\ell} \cdot
a(X_\ell(t_{k-1,\ell})) + 
b(X_\ell(t_{k-1,\ell})) \cdot V_{k,\ell}
\end{align}
for $k=1,\dots,2^\ell$.

The multilevel approach relies on a coupling of $X_\ell$ with 
$\ell \geq 1$ to an Euler scheme $\tX_{\ell-1}$ with step-size 
$2^{-(\ell-1)}$. Hence we choose 
\begin{equation}\label{g7}
\tV_{\ell-1} = (\tV_{1,\ell-1}, \dots, \tV_{2^{\ell-1},\ell-1}) 
\end{equation}
with $d$-dimensional random vectors $\tV_{k,\ell-1}$ on the
probability space introduced above, and we define, as before,
$\tX_{\ell-1}(t_{0,\ell-1}) = x_0$ and
\begin{align}\label{euler2}
\tX_{\ell-1}(t_{k,\ell-1})
=
\tX_{\ell-1}(t_{k-1,\ell-1}) 
\begin{aligned}[t]
&+ 
2^{-(\ell-1)} \cdot a\bigl(\tX_{\ell-1}(t_{k-1,\ell-1})\bigr)\\
&+ b\bigl(\tX_{\ell-1}(t_{k-1,\ell-1})\bigr) 
\cdot \tV_{k,\ell-1}
\end{aligned}
\end{align}
for $k=1,\dots,2^{\ell-1}$. Of course, a natural coupling between
$X_\ell$ and $\tX_{\ell-1}$ is induced
by
\begin{equation}\label{g4}
\phantom{\qquad\quad k=1,\dots, 2^{\ell-1}.}
\tV_{k,\ell-1} = V_{2k,\ell} + V_{2k-1,\ell},
\qquad\quad k=1,\dots, 2^{\ell-1}.
\end{equation}

Actually, the multilevel approach is based on
a hierarchy $V_0, (V_1,\tV_{0}), \dots,$ $(V_L,\tV_{L-1})$ with
corresponding Euler schemes, and the following properties are
most convenient for its analysis:
\begin{itemize}
\item[(i)]
For every $\ell \in \{0,\dots,L\}$ the random vectors
$V_{1,\ell}, \dots, V_{2^\ell,\ell}$ are iid 
with iid real-valued components
$V_{k,\ell}^{(1)},\dots,V_{k,\ell}^{(d)}$.
\item[(ii)]
For every $\ell \in \{1,\dots,L\}$ the random vectors
$\tV_{\ell-1}$ and $V_{\ell-1}$ coincide in distribution.
\end{itemize}

In order to obtain processes with continuous paths we extend $X_\ell$ and 
$\tX_{\ell-1}$ onto $[0,1]$ by piecewise linear interpolation.

\subsection{The Classical Euler Scheme}\label{s2.1}

In the vast majority of papers, $V_{k,\ell}$ and $\tV_{k,\ell-1}$ are
chosen as Brownian increments, i.e.,
\[
V_{k,\ell} = W(t_{k,\ell}) - W(t_{k-1,\ell})
\]
and
\[
\tV_{k,\ell-1} = W(t_{k,\ell-1}) - W(t_{k-1,\ell-1}),
\]
so that we have \eqref{g4}, 
and (i) and (ii) are satisfied. Error bounds for the Euler 
scheme w.r.t.\ various error criteria and under different sets of 
assumptions concerning the drift and diffusion coefficients are well known
in this case. In order to simulate the corresponding distributions
a generator for random numbers from $[0,1]$ has to be available.

\subsection{Random Bit Euler Schemes}\label{s2.2}

In the present paper we study the random bit quadrature problem for
SDEs, i.e., we consider algorithms that are only
allowed to use random bits instead of random numbers, see also
\cite{BN17,6924076,GHMR17,GHMR19,bookbits,O16}.
This excludes the use of Brownian increments.

Heuristics and extensive tests for finite precision random bit
multilevel algorithms for field programmable gate arrays
(FPGAs) are presented in \cite{6924076,bookbits,O16}.
In \cite{GHMR17} the random bit quadrature problem is studied for 
Gaussian random fields $X$, and relations to random bit
approximation of Gaussian measures are exploited.

Motivated by the weak error analysis of the Euler scheme,
the multilevel construction in \cite{BN17} is based on
iid random vectors $V_{1,L}, \dots, V_{2^L,L}$, each of which has 
iid components $V_{k,L}^{(1)},\dots,V_{k,L}^{(d)}$ with 
\[
2^{L/2-1} \cdot V_{1,L}^{(1)} + 1/2 \sim B(1,1/2).
\]
Moreover, the coupling is defined by \eqref{g4}, and (ii) is 
assumed to hold. It follows that (i) is satisfied as well, and 
\begin{equation}\label{g8}
2^{L/2-1} \cdot V^{(1)}_{1,\ell} + 2^{L-\ell-1} \sim
B(2^{L-\ell},1/2).
\end{equation}
See \cite[Sec.~3]{BN17} for error bounds, and
\cite[Sec.~4.1]{BN17} for the discussion of
fast generation of random quantities in this context.

A different construction is presented and analyzed
in \cite{GHMR19}. Here the starting point
is the approximation of the standard normal distribution 
based on random bits. Let $\Phi^{-1}$ denote the inverse 
of the distribution function of $N(0,1)$,
and let $U$ be uniformly distributed on
\begin{equation}\label{g6}
D_q = \left\{\sum_{i=1}^q b_i \cdot 2^{-i} + 2^{-(q+1)} : 
b_i \in \{0,1\} \ \text{for} \ i=1,\ldots,q\right\},
\end{equation}
where $q \in \N$. Obviously, $q$ random bits suffice to simulate the 
distribution $\nu_q$ of $\Phi^{-1} \circ U$,
which serves as an approximation of $N(0,1)$.
Further properties of $\nu_q$, in particular, error
bounds and the weak asymptotic optimality among all approximations
based on $q$ random bits, have been established in
\cite[Sec.~2.2]{GHMR17}. In the construction from \cite{GHMR19}, 
(i) is assumed to hold with
\begin{equation}\label{g9}
2^{\ell/2} \cdot V^{(1)}_{1,\ell} \sim \nu_L,
\end{equation}
and the coupling is again defined by \eqref{g4}.
Consequently, the analogon to (i) also holds for
the random vectors $\tV_{1,\ell-1},\dots,\tV_{2^{\ell-1},\ell-1}$,
but property (ii) is not satisfied, which
introduces an additional bias term in the multilevel analysis.
See \cite{GHMR19} for error and cost bounds; in particular,
a variant of the corresponding multilevel Euler algorithm, which
also employs Bakhvalov's trick, is shown to be almost worst case optimal 
the class of all Lipschitz continuous functionals 
$f \colon C([0,1],\R^r) \to \R$ with Lipschitz constant at most one.
Observe that the number of random bits
that are needed to simulate the distribution
of $V_\ell$ with $\ell=0$ or the joint distribution of $V_\ell$ and
$\tV_{\ell-1}$ with $\ell \geq 1$ is of the order $d \cdot L \cdot
2^\ell$. 

\subsection{The Random Bit L\'evy-Ciesielski Euler Scheme}\label{s2.3}

In the sequel we present a new construction of a random bit Euler
scheme, which is based on the L\'evy-Ciesielski representation of the 
Brownian motion. Hereby we get matching distributions across
the levels in the sense of (ii), but the iid-property (i) is not
satisfied. The main advantage of the new construction, 
compared to the approaches from \cite{BN17,GHMR19}, is that it is
well suited as the building block for an adaptive multilevel algorithm.

Consider the sequence of 
Schauder functions $s^{(i,j)}$ with $i=0$ and $j=1$ or $i \in \N$ and 
$j = 1 \dots, 2^{i-1}$. These functions are given by
\[
\phantom{\quad\qquad t \in [0,1],}
s^{(i,j)}(t) = \int_0^t h^{(i,j)}(u) \, \mathrm{d} u, 
\quad\qquad t \in [0,1],
\]
with Haar wavelets $h^{(0,1)}=1$ and
\[
h^{(i,j)} = 2^{(i-1)/2} \cdot
\left( 1_{I^{(i,j)}} - 1_{J^{(i,j)}} \right)
\] 
for $i \in \N$ and $j=1,\dots,2^{i-1}$, where
\[
I^{(i,j)}
=
\left[
(j-1) / 2^{i-1},
(j-1/2) / 2^{i-1}
\right[
\]
and
\[
J^{(i,j)}
=
\left[
(j-1/2) / 2^{i-1},
j / 2^{i-1}
\right[.
\]
The L\'evy-Ciesielski representation states that
\begin{equation}\label{g3}
W_\ell = 
s^{(0,1)} \cdot Z^{(0,1)} + 
\sum_{i=1}^\ell \sum_{j=1}^{2^{i-1}} s^{(i,j)} \cdot Z^{(i,j)}
\end{equation}
with an independent sequence $Z^{(0,1)},\dots$ of 
$d$-dimensional standard normally distributed random vectors
converges to a $d$-dimensional Brownian motion as $\ell \to \infty$,
e.g., in mean square and almost surely w.r.t.\ the supremum-norm.
We add that 
\begin{equation}\label{g15}
\phantom{\qquad\quad k =0,\dots,2^\ell,}
W_n(t_{k,\ell}) = W_\ell(t_{k,\ell}),
\qquad\quad k =0,\dots,2^\ell,
\end{equation}
for $\ell,n \in \N_0$ with $\ell < n$.
In this sense $W_\ell$ already yields the values of the Brownian motion
at the discretization~\eqref{g1}.

In a random bit approximation that corresponds to $W_\ell$ the number 
of bits that are spent for the individual terms
should depend on $i$ and $\ell$,
but 
not on the shift parameter $j$.
We spend 
\[
q^{(i)}_\ell =
2 \cdot (\ell+1-i)
\]
random bits for the approximation of the distribution of each of
the components of $Z^{(i,j)}$.
This choice is motivated by \cite[Thm.~2]{GHMR17},
which determines the weak asymptotics for random bit approximation of 
a Brownian bridge with respect to the $L_2$-norm.

Accordingly, we consider an independent sequence 
$U^{(0,1)}_\ell,\dots,U^{(\ell,2^{\ell-1})}_\ell$ 
of $d$-dimensional random vectors,
with iid components that
are uniformly distributed on 
\[
D^{(i)}_\ell = D_{q^{(i)}_\ell}.
\]
To normalize the variances we put 
\[
\sigma^{(i)}_\ell = 
2^{-q^{(i)}_\ell/2} \cdot 
\left(\sum_{x \in D^{(i)}_\ell} (\Phi^{-1}(x))^2\right)^{1/2}.
\]
Replacing $Z^{(i,j)}$ by 
\begin{align*}
Y^{(i,j)}_\ell = 1/ \sigma^{(i)}_\ell \cdot 
\Phi^{-1} \circ U^{(i,j)}_\ell
\end{align*}
in \eqref{g3}, where $\Phi^{-1} \circ U^{(i,j)}_\ell$ denotes the 
application of $\Phi^{-1}$ to every component of $U^{(i,j)}_\ell$,
we obtain
a random bit counterpart to $W_\ell$. 

Next, we turn to the Brownian increments, and we put
\[
\Delta^{(i,j)}_{k,\ell} = s^{(i,j)} (t_{k,\ell}) - s^{(i,j)} (t_{k-1,\ell}).
\]
We use \eqref{g5}
with
\[
V_{k,\ell} =
\Delta^{(0,1)}_{k,\ell} \cdot Y^{(0,1)}_\ell +
\sum_{i=1}^\ell \sum_{j=1}^{2^{i-1}} 
\Delta^{(i,j)}_{k,\ell} \cdot Y^{(i,j)}_\ell
\]
to approximate, in distribution,
the Brownian increments corresponding to \eqref{g1}.

\begin{lemma}\label{l1}
The components of $V_{k,\ell}$ have mean zero and variance
$2^{-\ell}$.
\end{lemma}

The normalization is crucial in the definition of 
the random vectors $Y^{(i,j)}_\ell$.
In fact, without this normalization the variances of 
the Brownian increments are not even matched asymptotically,
and thus one can not expect the Euler scheme to convergence, 
in any reasonable sense, to the true solution of the SDE.

\begin{lemma}\label{l2}
The components of 
\[
V^\prime_{k,\ell} =
\Delta^{(0,1)}_{k,\ell} \cdot \Phi^{-1} \circ U^{(0,1)}_\ell +
\sum_{i=1}^\ell \sum_{j=1}^{2^{i-1}} 
\Delta^{(i,j)}_{k,\ell} \cdot \Phi^{-1} \circ U^{(i,j)}_\ell
\]
have mean zero and variance at most $0.9 \cdot 2^{-\ell}$ 
for $\ell \geq 1$.
\end{lemma}

See the Appendix for the proofs of Lemma \ref{l1} and Lemma \ref{l2}. 

Let $\ell \geq 1$. For the multilevel construction we 
have to couple $V_\ell$ in a suitable way to a
random vector $\tV_{\ell-1}$
that approximates, in distribution, the Brownian increments with 
step-size $2^{-(\ell-1)}$. To this end we introduce the rounding function
\begin{align*}
T_q \colon {[0,1[} \to D_q, 
\quad x \mapsto \frac{\lfloor{2^q x}\rfloor}{2^q} + 2^{-(q+1)},
\end{align*}
see \eqref{g6}, and we put
\[
T^{(i)}_{\ell-1} = T_{q^{(i)}_{\ell-1}}
\]
for $i=0,\dots,\ell-1$ to obtain
\begin{equation}\label{g2}
U^{(i,j)}_{\ell-1}
\stackrel{\mathrm{d}}{=} 
T^{(i)}_{\ell-1} \circ U^{(i,j)}_\ell
\end{equation}
for $j=1$ if $i=0$ and for $j=1,\dots,2^{i-1}$ if $i \geq 1$.
We define
\[
\tY^{(i,j)}_{\ell-1} = 1/ \sigma^{(i)}_{\ell-1} \cdot 
\Phi^{-1} \circ 
T^{(i)}_{\ell-1} \circ U^{(i,j)}_\ell,
\]
and we use \eqref{g7} with
\[
\tV_{k,\ell-1} =
\Delta^{(0,1)}_{k,\ell-1} \cdot \tY^{(0,1)}_{\ell-1} +
\sum_{i=1}^{\ell-1} \sum_{j=1}^{2^{i-1}} 
\Delta^{(i,j)}_{k,\ell-1} \cdot \tY^{(i,j)}_{\ell-1}.
\]
Observe that $q^{(i)}_\ell - q^{(i)}_{\ell-1} = 2$.
Hence $\tV_{\ell-1}$ is, roughly speaking, obtained from 
$V_{\ell}$ by ignoring the two least important bits
in all of the relevant terms.

We stress that neither $V_\ell$ nor $\tV_{\ell-1}$ has
independent components, except for the trivial cases $\ell=0$ or 
$\ell=1$, respectively, so that (i) is not satisfied. 
On the other hand, \eqref{g2} implies that we have matching
distributions in the sense of (ii).

The number of random bits that are needed to simulate the distribution
of $V_\ell$ with $\ell=0$ or the joint distribution of $V_\ell$ and
$\tV_{\ell-1}$ with $\ell \geq 1$ is given by
\[
d \cdot 
\left( q^{(0)}_\ell + \sum_{i=1}^\ell 2^{i-1} \cdot q^{(i)}_\ell \right)
=
d \cdot \left( 2^{\ell+2} - 2 \right),
\]
which is easily verified by induction,
cf.~\cite[Thm.~2]{GHMR17}. Furthermore,
the arithmetic cost to compute
$V_\ell$, together with $\tV_{\ell-1}$ if $\ell \geq 1$, is of
the order $d \cdot 2^\ell$,
see, e.g., \cite[Sec.~2.2]{MR2900834}.

Let us discuss two important differences between the two constructions 
from \cite{BN17,GHMR19}, which have been discussed in Section \ref{s2.2}, 
and the construction based on the L\'evy-Ciesielski representation.

The maximal level $L$ has to be known in advance for the former
two constructions, see \eqref{g8} and \eqref{g9}, 
while this is not the case for the latter construction.
Due to this difference the L\'evy-Ciesielski based construction 
is well suited as the building block for an adaptive multilevel
algorithm.

On the other hand, analytic results are only available for
the constructions from \cite{BN17,GHMR19},
since we have (i) and its analogon for $\tV_{\ell-1}$ only in these
two cases.

\section{Adaptive Algorithms and Experiments}\label{s3}

We consider an adaptive multilevel algorithm with
either one of the following building blocks:
\begin{enumerate}
\item
The random bit Euler schemes with $V_\ell$ and $\tV_\ell$ based on
the L\'evy-Ciesielski representation, see Section \ref{s2.3}.
Notation: $X^\bit_\ell$ and $\tX^\bit_{\ell-1}$.
\item
The classical Euler schemes with $V_\ell$ and $\tV_\ell$ based on
the Brownian increments, see Section \ref{s2.1}.
Notation: $X^\num_\ell$ and $\tX^\num_{\ell-1}$.
\end{enumerate}
The number of calls to the random number generator
as well as
the number of arithmetic operations to jointly simulate 
$X^\ast_\ell$ and $\tX^\ast_{\ell-1}$ with $\ell \geq 1$
is of the order $2^\ell$ for both variants. 

In both cases we use the 
adaptive algorithm $A^\ast_\eps$ from \cite{MLMC}.
Here $\eps>0$ is an accuracy demand and input to the algorithm,
and the maximal level as well as the replication numbers per level
are determined adaptively.

We present numerical results for three different scalar
SDEs, i.e.,
$r=d=1$, where the solutions $\E(f(X))$ are known analytically.
For a fixed SDE and a fixed functional $f$ we put
\[
\bias_\ell^\ast = 
\E\bigl(f(X^\ast_\ell) - f(X^\ast_{\ell-1})\bigr)
\]
and
\[
\va_\ell^\ast = 
\var\bigl(f(X^\ast_\ell) - f(X^\ast_{\ell-1})\bigr)
\]
for $\ell \geq 1$. As key quantities we consider
the root mean squared error 
\[
\err^\ast_\eps =
\left( \E \left( A^\ast_\eps(f) - \E(f(X)) \right)^2 \right)^{1/2}
\]
of $A^\ast_\eps$, applied to $f$ for the particular SDE, and
the corresponding cost
\[
\cost_\eps^\ast = \E \left( C^\ast_\varepsilon(f)\right),
\]
where $C^\ast_\varepsilon(f)$ denotes the number of calls 
of the random number generator for $\eps > 0$.
All of these quantities can be 
approximated by simple Monte Carlo algorithms, and the
corresponding results will be presented together with asymptotic
confidence intervals with confidence level $0.95$ in the
sequel. The number of Monte Carlo replications for the 
data points and confidence intervals involving root mean squared
errors varies between $2 \cdot 10^3$ and $2 \cdot 10^4$.

\subsection{Geometric Brownian Motion}

Here we consider the geometric Brownian motion $X$ that solves 
\[
\phantom{\qquad\quad t \in [0,1]}
\mathrm{d}X(t) = 1/50 \cdot X(t)\,\mathrm{d}t + 
1/5 \cdot X(t)\,\mathrm{d}W(t), 
\qquad\quad t \in [0,1],
\] 
with initial value $x_0=1$, as well as the path-dependent
functional given by
\[
f(x) = \max_{0 \leq t \leq 1} x(t).
\]
Since $X(t) = \exp(W(t)/5)$, we obtain 
\[
\E(f(X)) = \E(\exp (|W(1)|/5)) = (2/\pi)^{1/2} \cdot \int_0^\infty
\exp(y/5 - y^2/2) \, \mathrm{d} y = 1.1819\dots.
\]

\begin{figure}
\centering
{
\begin{subfigure}[c]{0.45\textwidth}
\definecolor{mycolor1}{rgb}{0.80000,0.00000,0.10000}%
\definecolor{mycolor2}{rgb}{0.10000,0.13000,0.70000}%
\begin{tikzpicture}

\begin{axis}[%
width=0.69\textwidth,
height=3in,
scale only axis,
xmin=0,
xmax=11,
xtick={ 1,  2,  3,  4,  5,  6,  7,  8,  9, 10, 11},
xlabel style={font=\color{white!15!black}},
xlabel={level $\ell$},
ymin=-9.4,
ymax=-5.3,
ylabel style={font=\color{white!15!black}},
ylabel={$\log_2(|\bias_\ell^\ast|)$},
axis background/.style={fill=white},
axis x line*=bottom,
axis y line*=left,
xmajorgrids,
ymajorgrids,
legend style={legend cell align=left, align=left, draw=white!15!black}
]

\addplot [color=mycolor2, dashdotted, line width=0.7pt, mark=*, mark size =0.5pt, mark options={solid, mycolor2}, forget plot]
  table[row sep=crcr]{%
1	-6.20822639520065	0.0396046023204404	0.0407225834768239\\
2	-6.0075598979142	0.022960320024052	0.0233316483981927\\
3	-6.23866016141655	0.0195770790264591	0.0198463947666694\\
4	-6.55104548106845	0.0178113866585967	0.01803403615334\\
5	-6.90172407173644	0.0165096394485662	0.0167007585262597\\
6	-7.29645179829846	0.0160462994848789	0.0162267829032308\\
7	-7.7211864302694	0.0156371581517964	0.0158085059048583\\
8	-8.16156322379951	0.0152996051933449	0.0154635963020961\\
9	-8.64974696128559	0.0153782002967215	0.015543889739563\\
10	-9.08479663731246	0.014960607976624	0.0151173751133094\\
};

\addplot [color=mycolor2, draw=none,line width=0.7pt, mark=none, mark options={solid, mycolor2}, mark size=4pt, forget plot]
 plot [error bars/.cd, y dir = both, y explicit]
 table[row sep=crcr, y error plus index=2, y error minus index=3]{%
1	-6.20822639520065	0.0396046023204404	0.0407225834768239\\
2	-6.0075598979142	0.022960320024052	0.0233316483981927\\
3	-6.23866016141655	0.0195770790264591	0.0198463947666694\\
4	-6.55104548106845	0.0178113866585967	0.01803403615334\\
5	-6.90172407173644	0.0165096394485662	0.0167007585262597\\
6	-7.29645179829846	0.0160462994848789	0.0162267829032308\\
7	-7.7211864302694	0.0156371581517964	0.0158085059048583\\
8	-8.16156322379951	0.0152996051933449	0.0154635963020961\\
9	-8.64974696128559	0.0153782002967215	0.015543889739563\\
10	-9.08479663731246	0.014960607976624	0.0151173751133094\\
};
\addlegendimage{mycolor2, dashdotted, line width=0.7pt, mark=*,mark size=1.25pt}
\addlegendentry{$X_\ell^\bit$}
\addplot [color=mycolor1, dotted, line width=0.9pt, mark=*, mark size =0.5pt, mark options={solid, mycolor1}, forget plot]
  table[row sep=crcr]{%
1	-5.73144755207392	0.0191133525735534	0.0193699769299034\\
2	-5.90248437170158	0.0173983480213611	0.017610729256881\\
3	-6.16600157528221	0.0162627079619861	0.0164481205672118\\
4	-6.50076198160489	0.0155289462674553	0.0156979178765297\\
5	-6.87835434943069	0.0149853979963224	0.0151426878317373\\
6	-7.32471423304507	0.0147157569147103	0.0148674085671496\\
7	-7.76592233062842	0.0144496523798043	0.0145958417086796\\
8	-8.2202094137259	0.0142235125219656	0.0143651394119235\\
9	-8.70154812078656	0.0141388637147859	0.0142788015780653\\
10	-9.17507595354355	0.0141425476987216	0.0142825588568254\\
};

\addplot [color=mycolor1, draw=none, line width=0.9pt, mark=none, mark options={solid, mycolor1}, mark size =4pt, forget plot]
 plot [error bars/.cd, y dir = both, y explicit]
 table[row sep=crcr, y error plus index=2, y error minus index=3]{%
1	-5.73144755207392	0.0191133525735534	0.0193699769299034\\
2	-5.90248437170158	0.0173983480213611	0.017610729256881\\
3	-6.16600157528221	0.0162627079619861	0.0164481205672118\\
4	-6.50076198160489	0.0155289462674553	0.0156979178765297\\
5	-6.87835434943069	0.0149853979963224	0.0151426878317373\\
6	-7.32471423304507	0.0147157569147103	0.0148674085671496\\
7	-7.76592233062842	0.0144496523798043	0.0145958417086796\\
8	-8.2202094137259	0.0142235125219656	0.0143651394119235\\
9	-8.70154812078656	0.0141388637147859	0.0142788015780653\\
10	-9.17507595354355	0.0141425476987216	0.0142825588568254\\
};
\addlegendimage{mycolor1, dotted, line width=0.9pt, mark=*,mark size=1.25pt}
\addlegendentry{$X_\ell^\num$}
\end{axis}
\end{tikzpicture}%
\end{subfigure}
\quad
\begin{subfigure}[c]{0.45\textwidth}
\definecolor{mycolor1}{rgb}{0.80000,0.00000,0.10000}%
\definecolor{mycolor2}{rgb}{0.10000,0.13000,0.70000}%
\begin{tikzpicture}

\begin{axis}[%
width=0.69\textwidth,
height=3in,
scale only axis,
xmin=0,
xmax=11,
xtick={ 1,  2,  3,  4,  5,  6,  7,  8,  9, 10, 11},
xlabel style={font=\color{white!15!black}},
xlabel={level $\ell$},
ymin=-17.5,
ymax=-7.5,
ylabel style={font=\color{white!15!black}},
ylabel={$\log_2(\va_\ell^\ast)$},
axis background/.style={fill=white},
axis x line*=bottom,
axis y line*=left,
xmajorgrids,
ymajorgrids,
legend style={legend cell align=left, align=left, draw=white!15!black}
]

\addplot [color=mycolor2, dashdotted, line width=0.7pt, mark=*, mark size =0.5pt, mark options={solid, mycolor2}, forget plot]
  table[row sep=crcr]{%
1	-8.08268034482196	0.0133931791531232	0.0135186799763929\\
2	-9.27110246318266	0.0180563106286264	0.0182851649416715\\
3	-10.1966489942485	0.021123230847321	0.0214371081259443\\
4	-11.0959205140908	0.0210269743104252	0.0213379763783514\\
5	-12.0175644717539	0.0183480882094802	0.0185844470514525\\
6	-12.8896201061818	0.0175007259284694	0.017715629435564\\
7	-13.8140238266393	0.0147592256054558	0.0149117791568791\\
8	-14.7580834117833	0.0150369896714455	0.0151953701383345\\
9	-15.7195876510579	0.0341713452905807	0.0350003963852643\\
10	-16.669541067048	0.015417250200624	0.0155837867476514\\
};

\addplot [color=mycolor2, draw=none,line width=0.7pt, mark=none, mark options={solid, mycolor2}, mark size=4pt, forget plot]
 plot [error bars/.cd, y dir = both, y explicit]
 table[row sep=crcr, y error plus index=2, y error minus index=3]{%
1	-8.08268034482196	0.0133931791531232	0.0135186799763929\\
2	-9.27110246318266	0.0180563106286264	0.0182851649416715\\
3	-10.1966489942485	0.021123230847321	0.0214371081259443\\
4	-11.0959205140908	0.0210269743104252	0.0213379763783514\\
5	-12.0175644717539	0.0183480882094802	0.0185844470514525\\
6	-12.8896201061818	0.0175007259284694	0.017715629435564\\
7	-13.8140238266393	0.0147592256054558	0.0149117791568791\\
8	-14.7580834117833	0.0150369896714455	0.0151953701383345\\
9	-15.7195876510579	0.0341713452905807	0.0350003963852643\\
10	-16.669541067048	0.015417250200624	0.0155837867476514\\
};
\addlegendimage{mycolor2, dashdotted, line width=0.7pt, mark=*,mark size=1.25pt}
\addlegendentry{$X_\ell^\bit$}
\addplot [color=mycolor1, dotted, line width=0.9pt, mark=*, mark size =0.5pt, mark options={solid, mycolor1}, forget plot]
  table[row sep=crcr]{%
1	-9.25185787709839	0.0281639599371477	0.028724735769126\\
2	-9.86691119364221	0.0265726507629438	0.0270712845619538\\
3	-10.5898490938708	0.0251389927397199	0.0255848202748137\\
4	-11.3933201521721	0.0237461685319094	0.0241435709268583\\
5	-12.2518544675049	0.0227393859215077	0.0231035456901143\\
6	-13.1972356397171	0.021698098558069	0.0220294269479755\\
7	-14.132572316561	0.0214533183683514	0.0217771574456673\\
8	-15.086886998454	0.0208543360505136	0.021160214991264\\
9	-16.066872408147	0.020890087065915	0.0211970233915046\\
10	-17.0131726721821	0.0204905701638438	0.0207857954547528\\
};

\addplot [color=mycolor1, draw=none, line width=0.9pt, mark=none, mark options={solid, mycolor1}, mark size =4pt,forget plot]
 plot [error bars/.cd, y dir = both, y explicit]
 table[row sep=crcr, y error plus index=2, y error minus index=3]{%
1	-9.25185787709839	0.0281639599371477	0.028724735769126\\
2	-9.86691119364221	0.0265726507629438	0.0270712845619538\\
3	-10.5898490938708	0.0251389927397199	0.0255848202748137\\
4	-11.3933201521721	0.0237461685319094	0.0241435709268583\\
5	-12.2518544675049	0.0227393859215077	0.0231035456901143\\
6	-13.1972356397171	0.021698098558069	0.0220294269479755\\
7	-14.132572316561	0.0214533183683514	0.0217771574456673\\
8	-15.086886998454	0.0208543360505136	0.021160214991264\\
9	-16.066872408147	0.020890087065915	0.0211970233915046\\
10	-17.0131726721821	0.0204905701638438	0.0207857954547528\\
};
\addlegendimage{mycolor1, dotted, line width=0.9pt, mark=*,mark size=1.25pt}
\addlegendentry{$X_\ell^\num$}
\end{axis}
\end{tikzpicture}%
\end{subfigure}
}
\caption{Maximum of a geometric Brownian motion: 
bias and variance vs.\ level}
\label{f1}
\end{figure}
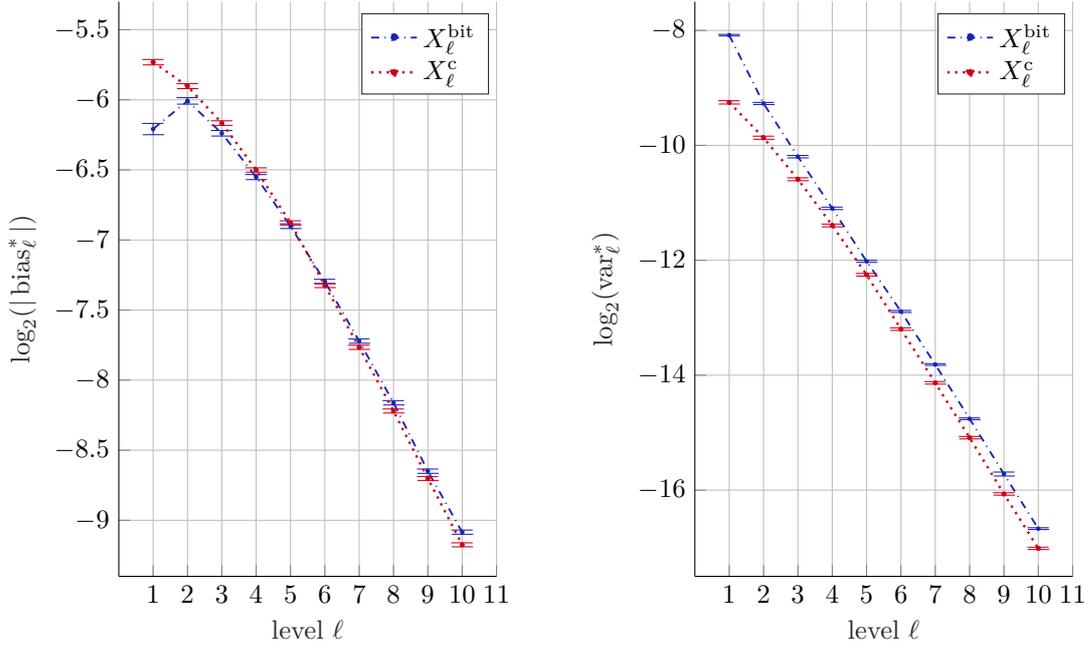

At first we compare the random bit Euler scheme $X^\bit_\ell$ and
the classical Euler scheme $X^\num_\ell$ in terms of their
bias and variance, see Figure \ref{f1}.
Since $f$ is Lipschitz continuous w.r.t.\ the supremum
norm on $C([0,1])$, we have the well-known upper bound
\[
\va_\ell^\num = O\left(\ell \cdot 2^{-\ell}\right)
\]
and, consequently,
\[
|\bias_\ell^\num| = O\left(\ell^{1/2} \cdot 2^{-\ell/2}\right).
\]
These upper bounds are very well reflected in the actual bias and 
variance decays, and we observe no essential difference between
the random bit and the classical Euler scheme.

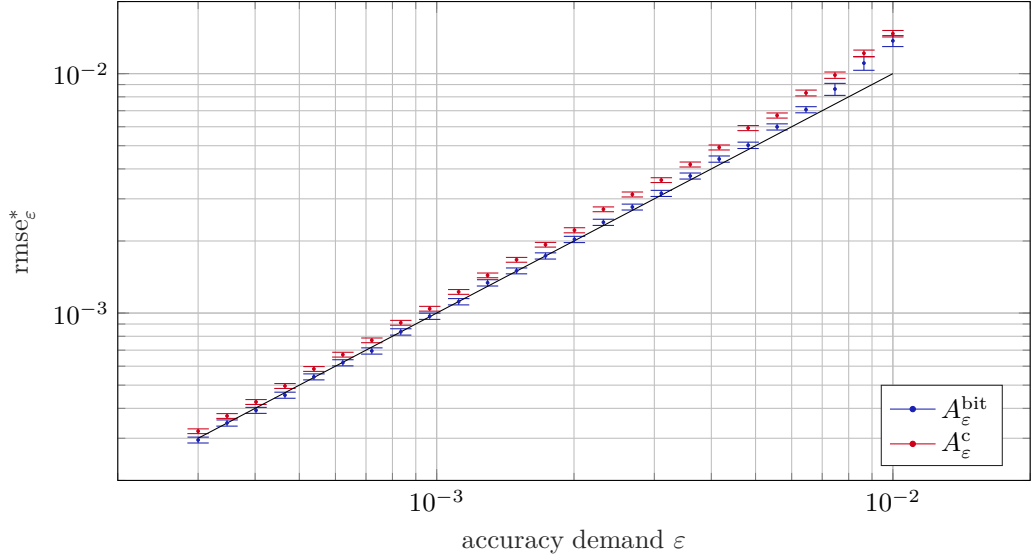
\begin{figure}
\centering
\definecolor{mycolor1}{rgb}{0.10000,0.13000,0.70000}%
\definecolor{mycolor2}{rgb}{0.80000,0.00000,0.10000}%
\begin{tikzpicture}

\begin{axis}[%
width=0.75\textwidth,
height=2.5in,
scale only axis,
xmode=log,
xmin=0.0002,
xmax=0.02,
xminorticks=true,
xlabel style={font=\color{white!15!black}},
xlabel={accuracy demand $\varepsilon$},
ymode=log,
ymin=0.0002,
ymax=0.02,
yminorticks=true,
ylabel style={font=\color{white!15!black}},
ylabel={$\err^\ast_\eps$},
axis background/.style={fill=white},
xmajorgrids,
xminorgrids,
ymajorgrids,
yminorgrids,
legend style={at={(0.97,0.03)}, anchor=south east, legend cell align=left, align=left, draw=white!15!black}
]
\addplot [color=mycolor1, draw=none, mark=*, mark size =0.5pt, mark options={solid, mycolor1}, forget plot]
  table[row sep=crcr]{%
0.01	0.0137092132101621\\
0.00864065606981118	0.0110641041444718\\
0.00746609373167647	0.00862394186046999\\
0.00645119481203894	0.00707240971375928\\
0.00557425556101787	0.00599405027582639\\
0.00481652251479877	0.00502299955637824\\
0.00416179145028782	0.00440004636134522\\
0.00359606085562177	0.00373867446190151\\
0.00310723250595386	0.00316328313150391\\
0.00268485274128848	0.00277363632988928\\
0.00231988891355635	0.00239515692272587\\
0.00200453622222083	0.00203164775257614\\
0.00173205080756888	0.00173358579726596\\
0.00149660553236414	0.00150130082850047\\
0.00129316536773352	0.00133820225272063\\
0.00111737971839762	0.00111603833486713\\
0.00096548938460563	0.000970651130624013\\
0.000834246171143089	0.0008347852254365\\
0.000720843424240427	0.000695045450580599\\
0.000622856010904651	0.000620339554931497\\
0.000538188457124165	0.000541752214145894\\
0.000465030135875223	0.000454210842134542\\
0.000401816546619536	0.000392620212118716\\
0.000347195858249866	0.000347599393431175\\
0.0003	0.000295103710106418\\
};
\addplot [color=mycolor1, draw=none, mark=none, mark size=4pt, mark options={solid, mycolor1}, forget plot]
 plot [error bars/.cd, y dir = both, y explicit]
 table[row sep=crcr, y error plus index=2, y error minus index=3]{%
0.01	0.0137092132101621	0.000699449041844729	0.000737108300358157\\
0.00864065606981118	0.0110641041444718	0.000691543405104973	0.000737751849030502\\
0.00746609373167647	0.00862394186046999	0.000482930160373884	0.000511628477623873\\
0.00645119481203894	0.00707240971375928	0.000205002568174225	0.000211124936305672\\
0.00557425556101787	0.00599405027582639	0.000174988360037204	0.000180252917172488\\
0.00481652251479877	0.00502299955637824	0.000151196080103593	0.000155890701815364\\
0.00416179145028782	0.00440004636134522	0.000129889750098506	0.000133842569891663\\
0.00359606085562177	0.00373867446190151	0.000106942332754904	0.000110092792322387\\
0.00310723250595386	0.00316328313150391	9.14246386586434e-05	9.41468245019114e-05\\
0.00268485274128848	0.00277363632988928	7.76823798581492e-05	7.99216834562904e-05\\
0.00231988891355635	0.00239515692272587	7.02089314306379e-05	7.23300729321818e-05\\
0.00200453622222083	0.00203164775257614	6.00134718823412e-05	6.18410353593951e-05\\
0.00173205080756888	0.00173358579726596	5.08719393711257e-05	5.24106073703451e-05\\
0.00149660553236414	0.00150130082850047	4.1818301663559e-05	4.30170065254044e-05\\
0.00129316536773352	0.00133820225272063	3.98174698321944e-05	4.10391239551785e-05\\
0.00111737971839762	0.00111603833486713	3.32525875430393e-05	3.4274263800072e-05\\
0.00096548938460563	0.000970651130624013	2.88170812418818e-05	2.96992039917645e-05\\
0.000834246171143089	0.0008347852254365	2.50740258625204e-05	2.58508565476361e-05\\
0.000720843424240427	0.000695045450580599	2.05007739019096e-05	2.1124121808165e-05\\
0.000622856010904651	0.000620339554931497	1.80050501159304e-05	1.85434997499602e-05\\
0.000538188457124165	0.000541752214145894	1.53943405144292e-05	1.58447701083631e-05\\
0.000465030135875223	0.000454210842134542	1.30799498206652e-05	1.34679533898511e-05\\
0.000401816546619536	0.000392620212118716	1.11609441200554e-05	1.14876370286668e-05\\
0.000347195858249866	0.000347599393431175	1.01740803119802e-05	1.04809897312326e-05\\
0.0003	0.000295103710106418	8.23179036927457e-06	8.46809894076631e-06\\
};
\addlegendimage{mycolor1,mark=*,mark size=1pt}
\addlegendentry{$A^\bit_\eps$}

\addplot [color=mycolor2, draw=none, mark=*, mark size =0.5pt, mark options={solid, mycolor2}, forget plot]
  table[row sep=crcr]{%
0.01	0.0146924170423824\\
0.00864065606981118	0.012168925161379\\
0.00746609373167647	0.0098602209982452\\
0.00645119481203894	0.00830747862104492\\
0.00557425556101787	0.00668819501641535\\
0.00481652251479877	0.00592635654594576\\
0.00416179145028782	0.00491890506633309\\
0.00359606085562177	0.0041751336519984\\
0.00310723250595386	0.00359269919437588\\
0.00268485274128848	0.00312982632365357\\
0.00231988891355635	0.00271362755789265\\
0.00200453622222083	0.00221832903450149\\
0.00173205080756888	0.00193046737285265\\
0.00149660553236414	0.00166960517493755\\
0.00129316536773352	0.00143754579173749\\
0.00111737971839762	0.00122591055569458\\
0.00096548938460563	0.00104298254587652\\
0.000834246171143089	0.000911230614879955\\
0.000720843424240427	0.000769719888073628\\
0.000622856010904651	0.000670643077091548\\
0.000538188457124165	0.000584395513067436\\
0.000465030135875223	0.00049617141704402\\
0.000401816546619536	0.000425391871804552\\
0.000347195858249866	0.000371808092239089\\
0.0003	0.000321108194894877\\
};

\addplot [color=mycolor2, draw=none, mark=none, mark size=4pt, mark options={solid, mycolor2}, forget plot]
 plot [error bars/.cd, y dir = both, y explicit]
 table[row sep=crcr, y error plus index=2, y error minus index=3]{%
0.01	0.0146924170423824	0.000451496366645222	0.000465817879038722\\
0.00864065606981118	0.012168925161379	0.000373096493494499	0.000384903268359195\\
0.00746609373167647	0.0098602209982452	0.000295870974556477	0.000305028058881576\\
0.00645119481203894	0.00830747862104492	0.000229508625234964	0.000236031982734572\\
0.00557425556101787	0.00668819501641535	0.000162396637006321	0.000166439182801553\\
0.00481652251479877	0.00592635654594576	0.00014225653386283	0.000145756308347209\\
0.00416179145028782	0.00491890506633309	0.000116849891239926	0.000119694078006626\\
0.00359606085562177	0.0041751336519984	9.87421064422567e-05	0.000101134630703439\\
0.00310723250595386	0.00359269919437588	8.24006130338596e-05	8.43354150458215e-05\\
0.00268485274128848	0.00312982632365357	7.38468031688008e-05	7.56318097193615e-05\\
0.00231988891355635	0.00271362755789265	6.20095677525438e-05	6.34600928729506e-05\\
0.00200453622222083	0.00221832903450149	5.35087585187576e-05	5.48317609119164e-05\\
0.00173205080756888	0.00193046737285265	4.31623012740804e-05	4.41496730195494e-05\\
0.00149660553236414	0.00166960517493755	3.83423009175199e-05	3.92437732609802e-05\\
0.00129316536773352	0.00143754579173749	3.2195793444553e-05	3.29335749776791e-05\\
0.00111737971839762	0.00122591055569458	2.79334200619411e-05	2.85849251387697e-05\\
0.00096548938460563	0.00104298254587652	2.3753192011357e-05	2.43069117716588e-05\\
0.000834246171143089	0.000911230614879955	2.08644785528056e-05	2.13535425898689e-05\\
0.000720843424240427	0.000769719888073628	1.71941248542017e-05	1.75870883301411e-05\\
0.000622856010904651	0.000670643077091548	1.53395373903215e-05	1.56987082095345e-05\\
0.000538188457124165	0.000584395513067436	1.33885547686821e-05	1.37025661791756e-05\\
0.000465030135875223	0.00049617141704402	1.15605766352489e-05	1.18364371206031e-05\\
0.000401816546619536	0.000425391871804552	9.83753650797125e-06	1.00704886147483e-05\\
0.000347195858249866	0.000371808092239089	8.47918140751774e-06	8.67711805570617e-06\\
0.0003	0.000321108194894877	7.23798885121505e-06	7.40494455043101e-06\\
};
\addlegendimage{mycolor2,mark=*,mark size=1pt}
\addlegendentry{$A^\num_\eps$}

\addplot [color=black, forget plot]
  table[row sep=crcr]{%
0.0003	0.0003\\
0.01	0.01\\
};
\end{axis}
\end{tikzpicture}%
\caption{Maximum of a geometric Brownian motion: 
root mean squared error vs.\ accuracy demand} 
\label{f2}
\end{figure}

Next we compare the multilevel algorithms $A^\bit_\eps$ and
$A^\num_\eps$. At first we relate the root mean squared error 
$\err^\ast_\eps$ to the accuracy demand~$\eps$, see Figure
\ref{f2}, where we consider $25$ different values of $\eps$ in a 
reasonable range. For both algorithms the root mean squared error 
is almost equal to the accuracy demand.

\begin{figure}
\centering
\input{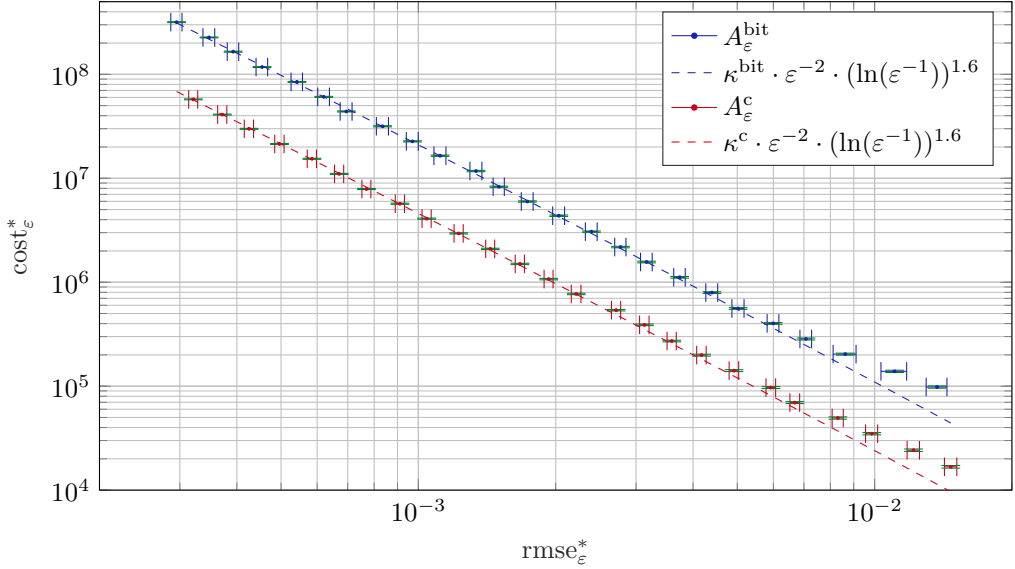}
\caption{Maximum of a geometric Brownian motion: 
cost vs.\ root mean squared error}
\label{f3}
\end{figure}

Finally we relate $\cost^\ast_\eps$ to the root mean squared error 
$\err_\eps^\ast$, see Figure~\ref{f3}, which is based on the same
data set as Figure~\ref{f2}.
We add that the confidence intervals
for $\cost^\ast_\eps$ in Figure~\ref{f3} are rather small and
hardly visible. Figure~\ref{f3} includes two graphs of
functions $\varepsilon \mapsto \kappa \cdot \varepsilon^{-2} \cdot (\ln
(\varepsilon^{-1}))^\gamma$ with parameters $\kappa > 0$ and $\gamma \in
\R$, which are fitted to the respective data by hand.
We obtain a log-exponent of $\gamma=1.6$ as a good fit for both algorithms.
The presence of a logarithmic term, i.e., $\gamma\neq 0$,
corresponds to the actual bias and variance decays.
The number of random bits is roughly $\kappa^\bit/\kappa^\num=4.57$ times 
larger than the number of random numbers for the same root mean squared 
error.  

\subsection{Ornstein-Uhlenbeck Process}\label{ornstein}

Here we consider the Ornstein-Uhlenbeck process $X$ that solves 
\[
\phantom{\qquad\quad t \in [0,1]}
\mathrm{d}X(t) = (2-X(t))\,\mathrm{d}t + \mathrm{d}W(t), 
\qquad\quad t \in [0,1],
\] 
with initial value $x_0=1$, as well as the path-independent
functional given by
\begin{equation}\label{g10}
f(x) = x(1).
\end{equation}
Since $X(t) = \exp(-t)+2(1-\exp(-t))+
\int_0^t \exp(-(t-s))\,\mathrm dW_s$, we obtain 
\[
\E(f(X)) = 2-\exp(-1)=1.6321\dots.
\]
As a major difference to the previous example we have
improved upper bounds for the classical Euler scheme $X^\num_\ell$,
namely,
\begin{equation}\label{eq100}
\va_\ell^\num = O\left(2^{-2\ell}\right)
\end{equation}
and, consequently,
\begin{equation}\label{eq101}
|\bias_\ell^\num| = O\left(2^{-\ell}\right).
\end{equation}

\begin{figure}
\centering
{
\begin{subfigure}[c]{0.45\textwidth}
\definecolor{mycolor1}{rgb}{0.80000,0.00000,0.10000}%
\definecolor{mycolor2}{rgb}{0.10000,0.13000,0.70000}%
\begin{tikzpicture}

\begin{axis}[%
width=0.69\textwidth,
height=3in,
scale only axis,
xmin=0,
xmax=11,
xtick={ 0,  1,  2,  3,  4,  5,  6,  7,  8,  9, 10, 11},
xlabel style={font=\color{white!15!black}},
xlabel={level $\ell$},
ymin=-13,
ymax=-1,
ylabel style={font=\color{white!15!black}},
ylabel={$\log_2(|\bias_\ell^\ast|)$},
axis background/.style={fill=white},
axis x line*=bottom,
axis y line*=left,
xmajorgrids,
ymajorgrids,
legend style={legend cell align=left, align=left, draw=white!15!black}
]

\addplot [color=mycolor2, dashdotted, line width=0.7pt, mark=*, mark size =0.5pt, mark options={solid, mycolor2}, forget plot]
  table[row sep=crcr]{%
1	-1.99757736370867	0.0161996535807802	0.0163836230480643\\
2	-3.90957328989027	0.0236802173777062	0.024075396976396\\
3	-5.19905278176599	0.0255731450212684	0.0260346462485392\\
4	-6.33802958484945	0.026362316945697	0.026853015075786\\
5	-7.3864380936022	0.0264090629771676	0.0269015191807283\\
6	-8.43756399202038	0.0271739399095363	0.0276956175359331\\
7	-9.43515962444531	0.0274978142347457	0.0280321239145849\\
8	-10.4350612191774	0.0275702948743231	0.0281074526091984\\
9	-11.4452529955533	0.0272109157129847	0.0277340277165123\\
10	-12.4360838558139	0.0271692903908107	0.0276907877936541\\
};

\addplot [color=mycolor2, draw=none,line width=0.7pt, mark=none, mark options={solid, mycolor2}, mark size=4pt, forget plot]
 plot [error bars/.cd, y dir = both, y explicit]
 table[row sep=crcr, y error plus index=2, y error minus index=3]{%
1	-1.99757736370867	0.0161996535807802	0.0163836230480643\\
2	-3.90957328989027	0.0236802173777062	0.024075396976396\\
3	-5.19905278176599	0.0255731450212684	0.0260346462485392\\
4	-6.33802958484945	0.026362316945697	0.026853015075786\\
5	-7.3864380936022	0.0264090629771676	0.0269015191807283\\
6	-8.43756399202038	0.0271739399095363	0.0276956175359331\\
7	-9.43515962444531	0.0274978142347457	0.0280321239145849\\
8	-10.4350612191774	0.0275702948743231	0.0281074526091984\\
9	-11.4452529955533	0.0272109157129847	0.0277340277165123\\
10	-12.4360838558139	0.0271692903908107	0.0276907877936541\\
};
\addlegendimage{mycolor2, dashdotted, line width=0.7pt, mark=*,mark size=1.25pt}
\addlegendentry{$X_\ell^\bit$}

\addplot [color=mycolor1, dotted, line width=0.9pt, mark=*, mark size =0.5pt, mark options={solid, mycolor1}, forget plot]
  table[row sep=crcr]{%
1	-2.01258822955798	0.0126824240143395	0.0127949019961733\\
2	-3.89675663468496	0.01781784794946	0.0180406600234679\\
3	-5.18752681282272	0.018763304709557	0.0190105544118628\\
4	-6.32821220659784	0.0192981842567086	0.0195598299239759\\
5	-7.3848803184355	0.0193549317116206	0.0196181289276174\\
6	-8.40914733350886	0.0193633097576758	0.0196267364341356\\
7	-9.43464938005103	0.0195416460095554	0.0198099810536174\\
8	-10.4379994417718	0.0194633354137785	0.0197295094512722\\
9	-11.4499909596033	0.0196360620002345	0.0199070142617863\\
10	-12.4279388611322	0.0193296660958744	0.0195921719431436\\
};

\addplot [color=mycolor1, draw=none, line width=0.9pt, mark=none, mark options={solid, mycolor1}, mark size =4pt, forget plot]
 plot [error bars/.cd, y dir = both, y explicit]
 table[row sep=crcr, y error plus index=2, y error minus index=3]{%
1	-2.01258822955798	0.0126824240143395	0.0127949019961733\\
2	-3.89675663468496	0.01781784794946	0.0180406600234679\\
3	-5.18752681282272	0.018763304709557	0.0190105544118628\\
4	-6.32821220659784	0.0192981842567086	0.0195598299239759\\
5	-7.3848803184355	0.0193549317116206	0.0196181289276174\\
6	-8.40914733350886	0.0193633097576758	0.0196267364341356\\
7	-9.43464938005103	0.0195416460095554	0.0198099810536174\\
8	-10.4379994417718	0.0194633354137785	0.0197295094512722\\
9	-11.4499909596033	0.0196360620002345	0.0199070142617863\\
10	-12.4279388611322	0.0193296660958744	0.0195921719431436\\
};
\addlegendimage{mycolor1, dotted, line width=0.9pt, mark=*,mark size=1.25pt}
\addlegendentry{$X_\ell^\num$}

\end{axis}
\end{tikzpicture}%
\end{subfigure}
\quad
\begin{subfigure}[c]{0.45\textwidth}
\definecolor{mycolor1}{rgb}{0.80000,0.00000,0.10000}%
\definecolor{mycolor2}{rgb}{0.10000,0.13000,0.70000}%
\begin{tikzpicture}

\begin{axis}[%
width=0.69\textwidth,
height=3in,
scale only axis,
xmin=0,
xmax=11,
xtick={ 0,  1,  2,  3,  4,  5,  6,  7,  8,  9, 10, 11},
xlabel style={font=\color{white!15!black}},
xlabel={level $\ell$},
ymin=-24,
ymax=-2,
ylabel style={font=\color{white!15!black}},
ylabel={$\log_2(\va_\ell^\ast)$},
axis background/.style={fill=white},
axis x line*=bottom,
axis y line*=left,
xmajorgrids,
ymajorgrids,
legend style={legend cell align=left, align=left, draw=white!15!black}
]

\addplot [color=mycolor2, dashdotted, line width=0.7pt, mark=*, mark size =0.5pt, mark options={solid, mycolor2}, forget plot]
  table[row sep=crcr]{%
1	-2.26427293018838	0.011765148511786	0.0118618824608236\\
2	-4.98532826248501	0.0130764290486596	0.0131960373046045\\
3	-7.34049439008516	0.0151917344956676	0.0153534090600225\\
4	-9.52996130462411	0.020492178199488	0.0207874501635938\\
5	-11.6216195501952	0.0285898609692374	0.0291679000037135\\
6	-13.6407228208382	0.0473716564715634	0.0489800932102735\\
7	-15.6014028079857	0.0929651202920461	0.0993707951045675\\
8	-17.5935377816138	0.182379066777983	0.208817063090002\\
9	-19.6521402307067	0.0665196530952059	0.0697356137190575\\
10	-21.6382609511393	0.0863621429369736	0.0918629606529855\\
};

\addplot [color=mycolor2, draw=none,line width=0.7pt, mark=none, mark options={solid, mycolor2}, mark size=4pt, forget plot]
 plot [error bars/.cd, y dir = both, y explicit]
 table[row sep=crcr, y error plus index=2, y error minus index=3]{%
1	-2.26427293018838	0.011765148511786	0.0118618824608236\\
2	-4.98532826248501	0.0130764290486596	0.0131960373046045\\
3	-7.34049439008516	0.0151917344956676	0.0153534090600225\\
4	-9.52996130462411	0.020492178199488	0.0207874501635938\\
5	-11.6216195501952	0.0285898609692374	0.0291679000037135\\
6	-13.6407228208382	0.0473716564715634	0.0489800932102735\\
7	-15.6014028079857	0.0929651202920461	0.0993707951045675\\
8	-17.5935377816138	0.182379066777983	0.208817063090002\\
9	-19.6521402307067	0.0665196530952059	0.0697356137190575\\
10	-21.6382609511393	0.0863621429369736	0.0918629606529855\\
};
\addlegendimage{mycolor2, dashdotted, line width=0.7pt, mark=*,mark size=1.25pt}
\addlegendentry{$X_\ell^\bit$}

\addplot [color=mycolor1, dotted, line width=0.9pt, mark=*, mark size =0.5pt, mark options={solid, mycolor1}, forget plot]
  table[row sep=crcr]{%
1	-3.00408265045563	0.0125497676263238	0.0126598946750036\\
2	-5.78628982721395	0.0125801572917403	0.0126908206935559\\
3	-8.21770103131	0.0125249584551153	0.0126346486163627\\
4	-10.417433431806	0.0126433377116175	0.0127551204040923\\
5	-12.5222405711518	0.0125962737705763	0.0127072221486433\\
6	-14.5695174943233	0.0126131192339489	0.0127243658735345\\
7	-16.5938900118611	0.0127235829804704	0.0128367954722606\\
8	-18.6122546768166	0.0125844271475337	0.0126951660141827\\
9	-20.6105712708618	0.0125636773109363	0.0126740496912312\\
10	-22.612151980834	0.0125264500120359	0.0126361664147296\\
};

\addplot [color=mycolor1, draw=none, line width=0.9pt, mark=none, mark options={solid, mycolor1}, mark size =4pt,forget plot]
 plot [error bars/.cd, y dir = both, y explicit]
 table[row sep=crcr, y error plus index=2, y error minus index=3]{%
1	-3.00408265045563	0.0125497676263238	0.0126598946750036\\
2	-5.78628982721395	0.0125801572917403	0.0126908206935559\\
3	-8.21770103131	0.0125249584551153	0.0126346486163627\\
4	-10.417433431806	0.0126433377116175	0.0127551204040923\\
5	-12.5222405711518	0.0125962737705763	0.0127072221486433\\
6	-14.5695174943233	0.0126131192339489	0.0127243658735345\\
7	-16.5938900118611	0.0127235829804704	0.0128367954722606\\
8	-18.6122546768166	0.0125844271475337	0.0126951660141827\\
9	-20.6105712708618	0.0125636773109363	0.0126740496912312\\
10	-22.612151980834	0.0125264500120359	0.0126361664147296\\
};
\addlegendimage{mycolor1, dotted, line width=0.9pt, mark=*,mark size=1.25pt}
\addlegendentry{$X_\ell^\num$}

\end{axis}
\end{tikzpicture}%
\end{subfigure}
}
\caption{Ornstein-Uhlenbeck process at final time: 
bias and variance vs.\ level}
\label{f4}
\end{figure}
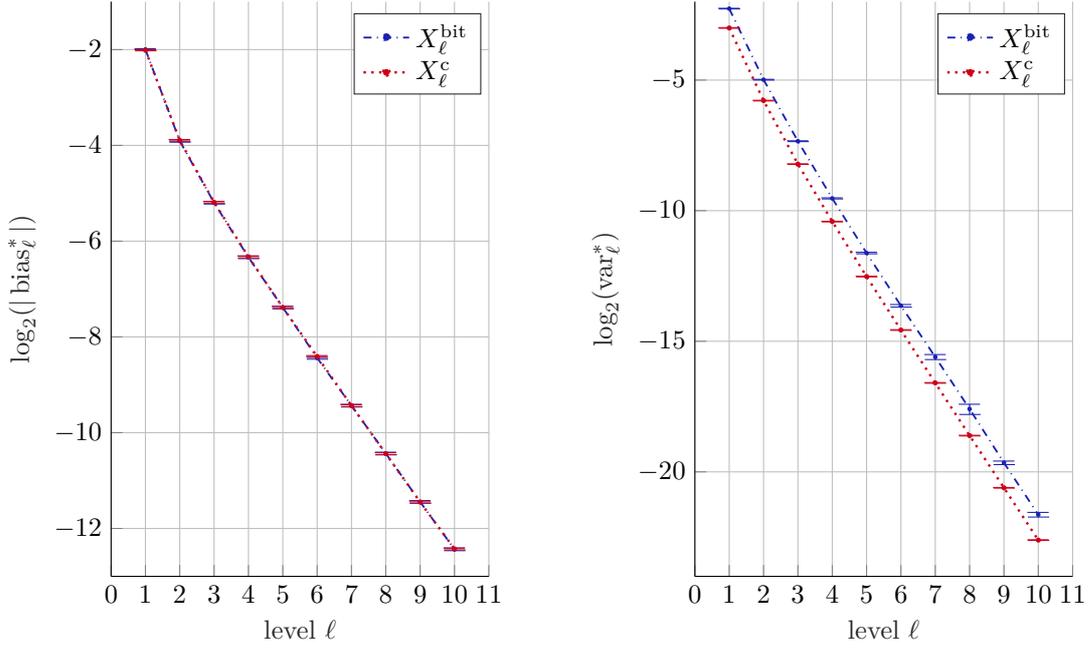

For the numerical experiments we proceed as in the previous section.
As for the geometric Brownian motion,
the upper bounds
\eqref{eq100} and \eqref{eq101}
are very well reflected in the actual bias and 
variance decays, and we observe no essential difference between
the random bit and the classical Euler scheme, see Figure~\ref{f4}.  

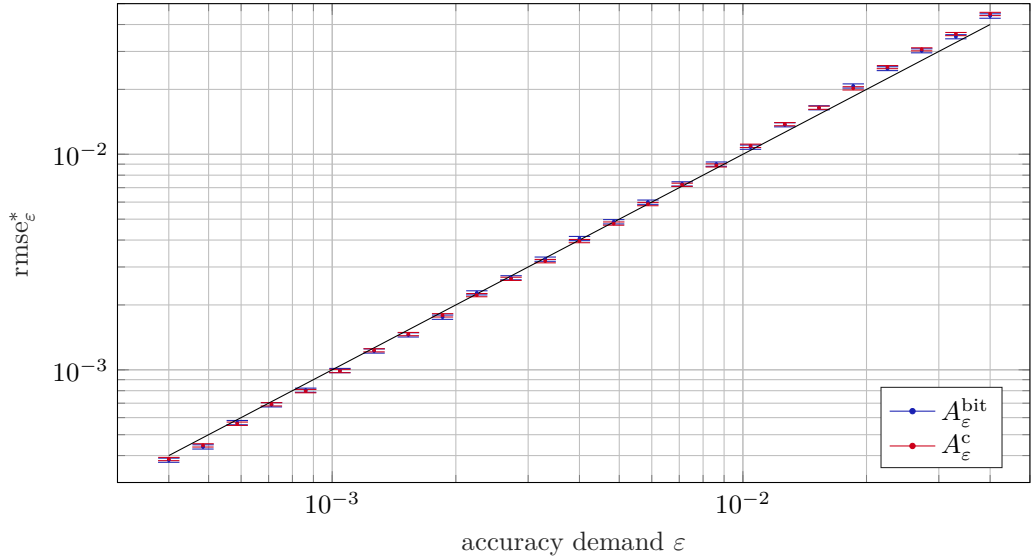
\begin{figure}
\centering
\definecolor{mycolor1}{rgb}{0.10000,0.13000,0.70000}%
\definecolor{mycolor2}{rgb}{0.80000,0.00000,0.10000}%
\begin{tikzpicture}

\begin{axis}[%
width=0.75\textwidth,
height=2.5in,
scale only axis,
xmode=log,
xmin=0.0003,
xmax=0.05,
xminorticks=true,
xlabel style={font=\color{white!15!black}},
xlabel={accuracy demand $\varepsilon$},
ymode=log,
ymin=0.0003,
ymax=0.05,
yminorticks=true,
ylabel style={font=\color{white!15!black}},
ylabel={$\err^\ast_\eps$},
axis background/.style={fill=white},
xmajorgrids,
xminorgrids,
ymajorgrids,
yminorgrids,
legend style={at={(0.97,0.03)}, anchor=south east, legend cell align=left, align=left, draw=white!15!black}
]
\addplot [color=mycolor1, draw=none, mark=*, mark size =0.5pt, mark options={solid, mycolor1}, forget plot]
  table[row sep=crcr]{%
0.04	0.0438216229188316\\
0.0330161674107207	0.0351221022390745\\
0.0272516827623185	0.0302166325293091\\
0.022493653007614	0.0249646988622227\\
0.0185663553344511	0.0207568137157996\\
0.0153247473982292	0.0164387249333607\\
0.0126491106406735	0.0137068836255226\\
0.0104406288627302	0.0107971342921387\\
0.00861773876012754	0.00900594001516391\\
0.00711311764015569	0.00729097070983763\\
0.00587119707048828	0.00599806565087823\\
0.00484611063451436	0.00486412767662031\\
0.004	0.00406673925010199\\
0.00330161674107207	0.00326308224728842\\
0.00272516827623184	0.00267460323633674\\
0.0022493653007614	0.00227638032604933\\
0.00185663553344511	0.00175533244817115\\
0.00153247473982291	0.00145291051879642\\
0.00126491106406735	0.00122310345874925\\
0.00104406288627301	0.000992625028469644\\
0.000861773876012754	0.00080485404996117\\
0.000711311764015569	0.000688545496766356\\
0.000587119707048828	0.00056912534450068\\
0.000484611063451436	0.000438981378915087\\
0.0004	0.000380870420363893\\
};
\addplot [color=mycolor1, draw=none, mark=none, mark size=4pt, mark options={solid, mycolor1}, forget plot]
 plot [error bars/.cd, y dir = both, y explicit]
 table[row sep=crcr, y error plus index=2, y error minus index=3]{%
0.04	0.0438216229188316	0.00102406312024209	0.0010485739958782\\
0.0330161674107207	0.0351221022390745	0.000770175580532821	0.000787447379820012\\
0.0272516827623185	0.0302166325293091	0.000640672909846055	0.000654554390969651\\
0.022493653007614	0.0249646988622227	0.000510375830507358	0.000521029989429662\\
0.0185663553344511	0.0207568137157996	0.000457567521965588	0.000467884223005529\\
0.0153247473982292	0.0164387249333607	0.000353836770451473	0.000361622386786886\\
0.0126491106406735	0.0137068836255226	0.0002956030214234	0.000302120100963948\\
0.0104406288627302	0.0107971342921387	0.000242520609594267	0.000248094644113947\\
0.00861773876012754	0.00900594001516391	0.000205342442963392	0.000210134958771556\\
0.00711311764015569	0.00729097070983763	0.000163994586003625	0.000167769166906039\\
0.00587119707048828	0.00599806565087823	0.000134853828990951	0.000137956286887511\\
0.00484611063451436	0.00486412767662031	0.000111140573381743	0.000113740118804748\\
0.004	0.00406673925010199	9.12375945268107e-05	9.33320451743458e-05\\
0.00330161674107207	0.00326308224728842	7.33019966917552e-05	7.49869410501445e-05\\
0.00272516827623184	0.00267460323633674	5.7926200470172e-05	5.92088453070377e-05\\
0.0022493653007614	0.00227638032604933	4.95584878959223e-05	5.06616979581995e-05\\
0.00185663553344511	0.00175533244817115	3.95813556504307e-05	4.04947170561064e-05\\
0.00153247473982291	0.00145291051879642	3.24946095727253e-05	3.32381777707847e-05\\
0.00126491106406735	0.00122310345874925	2.71659322876327e-05	2.77831705231114e-05\\
0.00104406288627301	0.000992625028469644	2.24240555644166e-05	2.29424766567645e-05\\
0.000861773876012754	0.00080485404996117	1.74435499964124e-05	1.78300728284559e-05\\
0.000711311764015569	0.000688545496766356	1.55167172463393e-05	1.58745511815199e-05\\
0.000587119707048828	0.00056912534450068	1.29796409750439e-05	1.32826496756579e-05\\
0.000484611063451436	0.000438981378915087	9.83008654847152e-06	1.00553123991484e-05\\
0.0004	0.000380870420363893	8.35651493428434e-06	8.54402181723398e-06\\
};
\addlegendimage{mycolor1,mark=*,mark size=1pt}
\addlegendentry{$A^\bit_\eps$}

\addplot [color=mycolor2, draw=none, mark=*, mark size =0.5pt, mark options={solid, mycolor2}, forget plot]
  table[row sep=crcr]{%
0.04	0.0447871758025057\\
0.0330161674107207	0.0361373273679131\\
0.0272516827623185	0.0307005448971403\\
0.022493653007614	0.0253815023789716\\
0.0185663553344511	0.0202333672277397\\
0.0153247473982292	0.0164472030438166\\
0.0126491106406735	0.0138022557769255\\
0.0104406288627302	0.0109506092462762\\
0.00861773876012754	0.00888010249677428\\
0.00711311764015569	0.00722047259317026\\
0.00587119707048828	0.00587770756308913\\
0.00484611063451436	0.00477140035966944\\
0.004	0.00395693934291957\\
0.00330161674107207	0.00319472445153283\\
0.00272516827623184	0.00264396606477777\\
0.0022493653007614	0.00222122199787137\\
0.00185663553344511	0.00179038584512651\\
0.00153247473982291	0.00146527053451028\\
0.00126491106406735	0.00123212428872559\\
0.00104406288627301	0.000987567215404683\\
0.000861773876012754	0.000796779230552105\\
0.000711311764015569	0.000692092972851261\\
0.000587119707048828	0.000562172513914238\\
0.000484611063451436	0.000446230977116488\\
0.0004	0.00038569692085697\\
};

\addplot [color=mycolor2, draw=none, mark=none, mark size=4pt, mark options={solid, mycolor2}, forget plot]
 plot [error bars/.cd, y dir = both, y explicit]
 table[row sep=crcr, y error plus index=2, y error minus index=3]{%
0.04	0.0447871758025057	0.000772880589442561	0.000786454282380618\\
0.0330161674107207	0.0361373273679131	0.000611311273837206	0.00062183192480407\\
0.0272516827623185	0.0307005448971403	0.000496858526738031	0.000505033085441541\\
0.022493653007614	0.0253815023789716	0.000420051921415063	0.000427121566944213\\
0.0185663553344511	0.0202333672277397	0.000350510060906691	0.000356690078097986\\
0.0153247473982292	0.0164472030438166	0.000276985326373921	0.000281730601371741\\
0.0126491106406735	0.0138022557769255	0.000233825462873118	0.000237855587712563\\
0.0104406288627302	0.0109506092462762	0.000194501074927715	0.000198018783627043\\
0.00861773876012754	0.00888010249677428	0.000147425433900214	0.000149914631251254\\
0.00711311764015569	0.00722047259317026	0.000123509042603213	0.000125658806264523\\
0.00587119707048828	0.00587770756308913	0.000103537380536389	0.00010539422104467\\
0.00484611063451436	0.00477140035966944	8.1661228367412e-05	8.30833900919174e-05\\
0.004	0.00395693934291957	6.80869383021506e-05	6.92792031228868e-05\\
0.00330161674107207	0.00319472445153283	5.77036482090482e-05	5.87652524275677e-05\\
0.00272516827623184	0.00264396606477777	4.52770209342706e-05	4.60660033160367e-05\\
0.0022493653007614	0.00222122199787137	3.77723954913959e-05	3.84259335513538e-05\\
0.00185663553344511	0.00179038584512651	3.18438427418053e-05	3.24205685894342e-05\\
0.00153247473982291	0.00146527053451028	2.52231280240944e-05	2.56649911224568e-05\\
0.00126491106406735	0.00123212428872559	2.05333925088521e-05	2.08814314176789e-05\\
0.00104406288627301	0.000987567215404683	1.72863012730254e-05	1.7594318925219e-05\\
0.000861773876012754	0.000796779230552105	1.40901762994869e-05	1.43438725125054e-05\\
0.000711311764015569	0.000692092972851261	1.21018232676124e-05	1.2317233905474e-05\\
0.000587119707048828	0.000562172513914238	1.00094986623684e-05	1.01909786387432e-05\\
0.000484611063451436	0.000446230977116488	7.80374974209146e-06	7.94267396076105e-06\\
0.0004	0.00038569692085697	6.6013483280454e-06	6.71631777610419e-06\\
};

\addlegendimage{mycolor2,mark=*,mark size=1pt}
\addlegendentry{$A^\num_\eps$}

\addplot [color=black, forget plot]
  table[row sep=crcr]{%
0.0004	0.0004\\
0.04	0.04\\
};
\end{axis}
\end{tikzpicture}%
\caption{Ornstein-Uhlenbeck process at final time: 
root mean squared error vs.\ accuracy demand} 
\label{f5}
\end{figure}

For both multilevel algorithms the root mean squared error is again
almost equal to the accuracy demand, see Figure~\ref{f5}.

\begin{figure}
\centering
\input{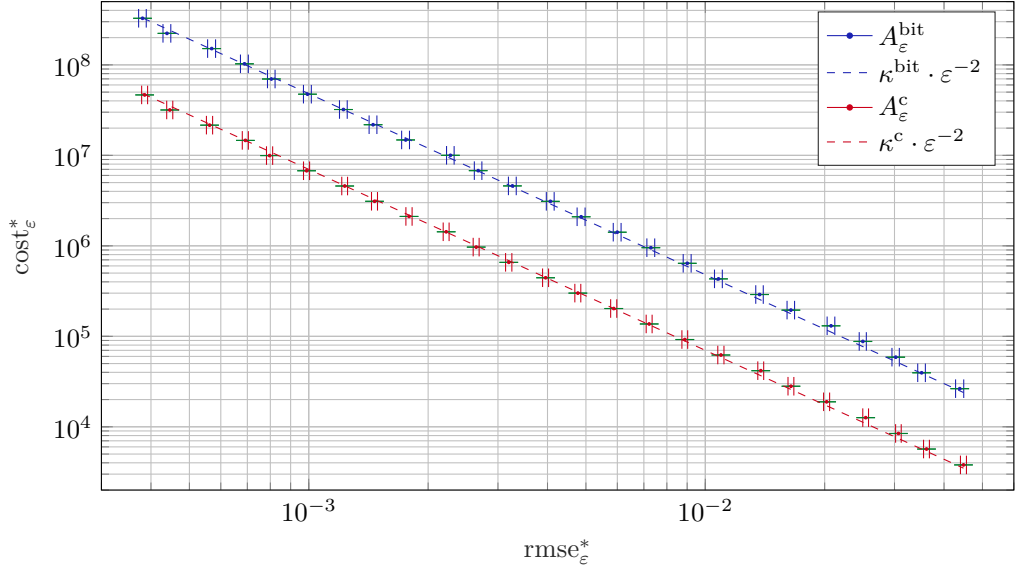}
\caption{Ornstein-Uhlenbeck process at final time: 
cost vs.\ root mean squared error}
\label{f6}
\end{figure}

Due to the improved upper bounds for the variance and bias
it is natural to expect that 
$\cost^\ast_\eps$ is proportional to
${(\err_\eps^\ast)}^{-2}$.
This is in line with the numerical results in Figure~\ref{f6}.
Furthermore, we have $\kappa^\bit/\kappa^\num=6.85$.

\subsection{Cox-Ingersoll-Ross Process}

Here we consider the Cox-Ingersoll-Ross process $X$ that solves 
\[
\phantom{\qquad\quad t \in [0,1]}
\mathrm{d}X(t) = (3/2-X(t))\,\mathrm{d}t + 2 \cdot \sqrt{X(t)}\,
\mathrm{d}W(t), 
\qquad\quad t \in [0,1],
\]
with initial value $x_0=1$, as well as $f$ given by \eqref{g10}.
We have
\[
\E(f(X)) = \exp(-1)+\frac{3}{2}\left(1-\exp(-1)\right)
=1.3160\dots,
\]
see, e.g., \cite[Eqn.~(19)]{MR785475}.
To get a well-defined variant of the Euler scheme we take the positive 
part in every Euler step, i.e., we take the maximum with $0$ of the 
right-hand side in \eqref{euler1} and \eqref{euler2}.

Furthermore, we compare this Euler scheme with a truncated Milstein 
scheme, which is proposed and analyzed in \cite{MR3732573}.
For this scheme the right-hand side of \eqref{euler1},
and similarly also for \eqref{euler2},
is replaced by
$\Theta_\hh(X_\ell(t_{k-1,\ell}), 2^{-\ell}, V_{k,\ell})$,
where
\begin{equation*}
\Theta_\hh(x,h,w)
=\max\left(0,\left(\max
\left(\sqrt{h},\sqrt{\max(h,x)}+w\right)\right)^2+(1/2-x) \cdot h\right).
\end{equation*}
The resulting schemes are denoted by
$X_\ell^{\num,\hh}$
and
$X_\ell^{\bit,\hh}$.
For the Euler scheme no polynomial strong convergence rate is 
known. 
For the truncated Milstein scheme the strong convergence result
from \cite[Thm.~1]{MR3732573} implies
\begin{equation}\label{varcox}
\va_\ell^{\num,\hh} = O\left(2^{-\ell/2+\varepsilon\cdot \ell}\right)
\end{equation}
and
\begin{equation}\label{biacox}
|\bias_\ell^{\num,\hh}| = O\left(2^{-\ell/2+\varepsilon\cdot \ell}\right)
\end{equation}
for every $\varepsilon>0$.
This strong convergence rate is the best known convergence rate 
for the Cox-Ingersoll-Ross process, see \cite[Fig.~1.1]{MR3744680}.

For the numerical experiments we proceed as in the previous sections.
The decay of the bias for all four variants is similar to the decay of 
the bias for the Ornstein-Uhlenbeck process.
The decay of the variance for both variants based on the Euler 
scheme is similar to the decay of the variance for the 
geometric Brownian motion. The decay of the variance for both versions 
based on the truncated Milstein scheme is similar and substantially
faster. Note that the upper bounds \eqref{varcox} and \eqref{biacox} 
seem to be too pessimistic, cf.\ the conjecture in
\cite[Fig.~5]{MR3732573}.

\begin{figure}
\centering
{
\begin{subfigure}[c]{0.45\textwidth}
\definecolor{mycolor1}{rgb}{0.80000,0.00000,0.10000}%
\definecolor{mycolor2}{rgb}{0.10000,0.13000,0.70000}%
\definecolor{mycolor4}{rgb}{0.50000,0.00000,0.80000}%
\definecolor{mycolor3}{rgb}{0.85000,0.45000,0.00000}%
\begin{tikzpicture}

\begin{axis}[%
width=0.69\textwidth,
height=3in,
scale only axis,
xmin=0,
xmax=11,
xtick={ 0,  1,  2,  3,  4,  5,  6,  7,  8,  9, 10, 11},
xlabel style={font=\color{white!15!black}},
xlabel={level $\ell$},
ymin=-11,
ymax=-2,
ylabel style={font=\color{white!15!black}},
ylabel={$\log_2(|\bias_\ell^\ast|)$},
axis background/.style={fill=white},
axis x line*=bottom,
axis y line*=left,
xmajorgrids,
ymajorgrids,
legend style={legend cell align=left, align=left, draw=white!15!black}
]
\addplot [color=mycolor2, dashdotted, line width=0.7pt, mark=*, mark size =0.5pt, mark options={solid, mycolor2}, forget plot]
  table[row sep=crcr]{%
1	-2.07268316397498	0.0340331347994516	0.0348554120257862\\
2	-3.28054245259389	0.0503743268081802	0.0521970746826872\\
3	-4.00670362479821	0.0534103529214089	0.0554639400603261\\
4	-4.76571218533961	0.0592589273210367	0.0617976476844442\\
5	-5.63511774636743	0.0729238591786352	0.0768070956525158\\
6	-6.4212573834318	0.0866384885802693	0.0921757071048726\\
7	-7.18552941575832	0.102401466559479	0.110228977403908\\
8	-8.04632916080623	0.129929978127246	0.142800206660171\\
9	-8.79976269631762	0.154423674157965	0.172956449733777\\
10	-9.6215873787173	0.189523379739487	0.218240636698884\\
};
\addplot [color=mycolor2, draw=none,line width=0.7pt, mark=none, mark options={solid, mycolor2}, mark size=4pt, forget plot]
 plot [error bars/.cd, y dir = both, y explicit]
 table[row sep=crcr, y error plus index=2, y error minus index=3]{%
1	-2.07268316397498	0.0340331347994516	0.0348554120257862\\
2	-3.28054245259389	0.0503743268081802	0.0521970746826872\\
3	-4.00670362479821	0.0534103529214089	0.0554639400603261\\
4	-4.76571218533961	0.0592589273210367	0.0617976476844442\\
5	-5.63511774636743	0.0729238591786352	0.0768070956525158\\
6	-6.4212573834318	0.0866384885802693	0.0921757071048726\\
7	-7.18552941575832	0.102401466559479	0.110228977403908\\
8	-8.04632916080623	0.129929978127246	0.142800206660171\\
9	-8.79976269631762	0.154423674157965	0.172956449733777\\
10	-9.6215873787173	0.189523379739487	0.218240636698884\\
};
\addlegendimage{mycolor2, dashdotted, line width=0.7pt, mark=*,mark size=1.25pt}
\addlegendentry{$X_\ell^\bit$}

\addplot [color=mycolor1, dotted, line width=0.9pt, mark=*, mark size =0.5pt, mark options={solid, mycolor1}, forget plot]
  table[row sep=crcr]{%
1	-2.5151660435873	0.036502257868567	0.0374498446189482\\
2	-3.21776099143895	0.047126210802209	0.0487177416149529\\
3	-3.81657464146885	0.0510969067525844	0.0529733013592435\\
4	-4.55101939333906	0.0604005840746931	0.0630402612106886\\
5	-5.48749952971499	0.082068246562387	0.087019816939951\\
6	-6.19283677863377	0.0951841508767961	0.101910471669118\\
7	-7.27799999447485	0.141095934441515	0.156406000471873\\
8	-8.24038533838145	0.191376041099831	0.220702106787048\\
9	-8.89921507345723	0.213390293946398	0.250525056111337\\
10	-9.49861453824612	0.226576098919947	0.268912011827105\\
};
\addplot [color=mycolor1, draw=none, line width=0.9pt, mark=none, mark options={solid, mycolor1}, mark size =4pt, forget plot]
 plot [error bars/.cd, y dir = both, y explicit]
 table[row sep=crcr, y error plus index=2, y error minus index=3]{%
1	-2.5151660435873	0.036502257868567	0.0374498446189482\\
2	-3.21776099143895	0.047126210802209	0.0487177416149529\\
3	-3.81657464146885	0.0510969067525844	0.0529733013592435\\
4	-4.55101939333906	0.0604005840746931	0.0630402612106886\\
5	-5.48749952971499	0.082068246562387	0.087019816939951\\
6	-6.19283677863377	0.0951841508767961	0.101910471669118\\
7	-7.27799999447485	0.141095934441515	0.156406000471873\\
8	-8.24038533838145	0.191376041099831	0.220702106787048\\
9	-8.89921507345723	0.213390293946398	0.250525056111337\\
10	-9.49861453824612	0.226576098919947	0.268912011827105\\
};
\addlegendimage{mycolor1, dotted, line width=0.9pt, mark=*,mark size=1.25pt}
\addlegendentry{$X_\ell^\num$}

\addplot [color=mycolor4, dashed, line width=0.6pt, mark=*, mark size =0.5pt, mark options={solid, mycolor4}, forget plot]
  table[row sep=crcr]{%
1	-2.12158257019083	0.0446735324496705	0.046101188832139\\
2	-2.48571211138526	0.0278920332572454	0.0284419265726408\\
3	-3.66083239557221	0.0311205445433629	0.0318066773151084\\
4	-5.01175969865958	0.0420627242208313	0.0433260148962198\\
5	-6.01821237941253	0.0449576478307705	0.0464038167384144\\
6	-6.9880020096293	0.0491456348167807	0.0508790098613465\\
7	-7.96553874370152	0.0543392808151344	0.0564663540675641\\
8	-8.84791615510144	0.0589319203942402	0.0614421016183595\\
9	-9.6799641858991	0.0624224973502852	0.0652460211373462\\
10	-10.4979152587159	0.0660895917629194	0.0692631032150288\\
};
\addplot [color=mycolor4, draw=none, line width=0.6pt, mark=none, mark options={solid, mycolor4}, mark size = 4pt, forget plot]
 plot [error bars/.cd, y dir = both, y explicit]
 table[row sep=crcr, y error plus index=2, y error minus index=3]{%
1	-2.12158257019083	0.0446735324496705	0.046101188832139\\
2	-2.48571211138526	0.0278920332572454	0.0284419265726408\\
3	-3.66083239557221	0.0311205445433629	0.0318066773151084\\
4	-5.01175969865958	0.0420627242208313	0.0433260148962198\\
5	-6.01821237941253	0.0449576478307705	0.0464038167384144\\
6	-6.9880020096293	0.0491456348167807	0.0508790098613465\\
7	-7.96553874370152	0.0543392808151344	0.0564663540675641\\
8	-8.84791615510144	0.0589319203942402	0.0614421016183595\\
9	-9.6799641858991	0.0624224973502852	0.0652460211373462\\
10	-10.4979152587159	0.0660895917629194	0.0692631032150288\\
};
\addlegendimage{mycolor4, dashed, line width=0.6pt, mark=*,mark size=1.25pt}
\addlegendentry{$X_\ell^{\bit,\hh}$}

\addplot [color=mycolor3, dashdotted, line width=0.6pt, mark=*, mark size =0.5pt, mark options={solid, mycolor3}, forget plot]
  table[row sep=crcr]{%
1	-2.73200450604896	0.0518594586548109	0.0537933429252875\\
2	-2.47578373105388	0.0216655835351682	0.0219959120869744\\
3	-3.61479043518457	0.0239278932375333	0.0243314532205439\\
4	-4.84751839916644	0.0293432469864987	0.0299524775233255\\
5	-5.90163822363867	0.0334105426504046	0.0342026585397512\\
6	-6.82490760164295	0.035499620933976	0.0363952249881603\\
7	-7.75667720861256	0.0390089487891432	0.0400930914842235\\
8	-8.62388633080856	0.0415101328562795	0.0427399613716943\\
9	-9.42215807694239	0.0432059866579237	0.0445399774836854\\
10	-10.1935189200055	0.0436419690516026	0.0450034441446316\\
};
\addplot [color=mycolor3, draw=none, line width=0.6pt, mark=none, mark options={solid, mycolor3}, mark size = 4pt, forget plot]
 plot [error bars/.cd, y dir = both, y explicit]
 table[row sep=crcr, y error plus index=2, y error minus index=3]{%
1	-2.73200450604896	0.0518594586548109	0.0537933429252875\\
2	-2.47578373105388	0.0216655835351682	0.0219959120869744\\
3	-3.61479043518457	0.0239278932375333	0.0243314532205439\\
4	-4.84751839916644	0.0293432469864987	0.0299524775233255\\
5	-5.90163822363867	0.0334105426504046	0.0342026585397512\\
6	-6.82490760164295	0.035499620933976	0.0363952249881603\\
7	-7.75667720861256	0.0390089487891432	0.0400930914842235\\
8	-8.62388633080856	0.0415101328562795	0.0427399613716943\\
9	-9.42215807694239	0.0432059866579237	0.0445399774836854\\
10	-10.1935189200055	0.0436419690516026	0.0450034441446316\\
};
\addlegendimage{mycolor3, dashdotted, line width=0.6pt, mark=*,mark size=1.25pt}
\addlegendentry{$X_\ell^{\num,\hh}$}

\end{axis}
\end{tikzpicture}%
\end{subfigure}
\quad
\begin{subfigure}[c]{0.45\textwidth}
\definecolor{mycolor1}{rgb}{0.80000,0.00000,0.10000}%
\definecolor{mycolor2}{rgb}{0.10000,0.13000,0.70000}%
\definecolor{mycolor4}{rgb}{0.50000,0.00000,0.80000}%
\definecolor{mycolor3}{rgb}{0.85000,0.45000,0.00000}%
\begin{tikzpicture}

\begin{axis}[%
width=0.69\textwidth,
height=3in,
scale only axis,
xmin=0,
xmax=11,
xtick={ 0,  1,  2,  3,  4,  5,  6,  7,  8,  9, 10, 11},
xlabel style={font=\color{white!15!black}},
xlabel={level $\ell$},
ymin=-17,
ymax=2,
ylabel style={font=\color{white!15!black}},
ylabel={$\log_2(\va_\ell^\ast)$},
axis background/.style={fill=white},
axis x line*=bottom,
axis y line*=left,
xmajorgrids,
ymajorgrids,
legend style={legend cell align=left, align=left, draw=white!15!black}
]

\addplot [color=mycolor2, dashdotted, line width=0.7pt, mark=*, mark size =0.5pt, mark options={solid, mycolor2}, forget plot]
  table[row sep=crcr]{%
1	-0.254645140792273	0.0136783606583421	0.01380928913622\\
2	-1.52244563559472	0.0167156696285098	0.0169116169306345\\
3	-2.80285222718194	0.0222419373545577	0.0225902162832408\\
4	-4.01515652352164	0.0261664979663445	0.026649866331061\\
5	-5.14148158270766	0.02489856039817	0.0253358264312391\\
6	-6.20268550374997	0.0258596646761928	0.0263316608111097\\
7	-7.2329816869276	0.0205919440255382	0.0208901190091311\\
8	-8.23968735738319	0.0201874335217092	0.0204739271121692\\
9	-9.22333927365911	0.0194895979059915	0.019756495678525\\
10	-10.2402279939245	0.0198236986737097	0.0200998904656924\\
};
\addplot [color=mycolor2, draw=none,line width=0.7pt, mark=none, mark options={solid, mycolor2}, mark size=4pt, forget plot]
 plot [error bars/.cd, y dir = both, y explicit]
 table[row sep=crcr, y error plus index=2, y error minus index=3]{%
1	-0.254645140792273	0.0136783606583421	0.01380928913622\\
2	-1.52244563559472	0.0167156696285098	0.0169116169306345\\
3	-2.80285222718194	0.0222419373545577	0.0225902162832408\\
4	-4.01515652352164	0.0261664979663445	0.026649866331061\\
5	-5.14148158270766	0.02489856039817	0.0253358264312391\\
6	-6.20268550374997	0.0258596646761928	0.0263316608111097\\
7	-7.2329816869276	0.0205919440255382	0.0208901190091311\\
8	-8.23968735738319	0.0201874335217092	0.0204739271121692\\
9	-9.22333927365911	0.0194895979059915	0.019756495678525\\
10	-10.2402279939245	0.0198236986737097	0.0200998904656924\\
};
\addlegendimage{mycolor2, dashdotted, line width=0.7pt, mark=*,mark size=1.25pt}
\addlegendentry{$X_\ell^\bit$}

\addplot [color=mycolor1, dotted, line width=0.9pt, mark=*, mark size =0.5pt, mark options={solid, mycolor1}, forget plot]
  table[row sep=crcr]{%
1	-0.935040394893406	0.0183643076889811	0.0186010872942936\\
2	-1.59246705831163	0.0230793807842846	0.0234546017433885\\
3	-2.55268859100003	0.0242705287953915	0.0246858296699428\\
4	-3.52956141408984	0.0240714060386491	0.0244798628534082\\
5	-4.49615246344172	0.0225400406033325	0.0228977931694407\\
6	-5.46578229697507	0.0222125074181116	0.0225598580741204\\
7	-6.45380412403448	0.0214009225538927	0.0217231698168607\\
8	-7.44786171248995	0.0199025642046067	0.0201809734087242\\
9	-8.42882373977768	0.0205237689533746	0.0208199585972366\\
10	-9.44109964672724	0.0191356936480851	0.0193929223234157\\
};
\addplot [color=mycolor1, draw=none, line width=0.9pt, mark=none, mark options={solid, mycolor1}, mark size =4pt, forget plot]
 plot [error bars/.cd, y dir = both, y explicit]
 table[row sep=crcr, y error plus index=2, y error minus index=3]{%
1	-0.935040394893406	0.0183643076889811	0.0186010872942936\\
2	-1.59246705831163	0.0230793807842846	0.0234546017433885\\
3	-2.55268859100003	0.0242705287953915	0.0246858296699428\\
4	-3.52956141408984	0.0240714060386491	0.0244798628534082\\
5	-4.49615246344172	0.0225400406033325	0.0228977931694407\\
6	-5.46578229697507	0.0222125074181116	0.0225598580741204\\
7	-6.45380412403448	0.0214009225538927	0.0217231698168607\\
8	-7.44786171248995	0.0199025642046067	0.0201809734087242\\
9	-8.42882373977768	0.0205237689533746	0.0208199585972366\\
10	-9.44109964672724	0.0191356936480851	0.0193929223234157\\
};
\addlegendimage{mycolor1, dotted, line width=0.9pt, mark=*,mark size=1.25pt}
\addlegendentry{$X_\ell^\num$}

\addplot [color=mycolor4, dashed, line width=0.6pt, mark=*, mark size =0.5pt, mark options={solid, mycolor4}, forget plot]
  table[row sep=crcr]{%
1	0.443205402663294	0.0157657391566289	0.0159399321312996\\
2	-1.66103997818893	0.0207128481597494	0.0210145605628551\\
3	-3.69201263493529	0.0310795472356711	0.0317638534657791\\
4	-5.51352883566145	0.0476674656814913	0.0492964000556411\\
5	-7.33147619048305	0.0628449140541454	0.0657076640802412\\
6	-9.00985125227136	0.103051617135609	0.110982737950412\\
7	-10.6698356471261	0.0574282865647753	0.0598094023026068\\
8	-12.1958601662628	0.0613471402384658	0.0640720803007842\\
9	-13.6904075895865	0.0819117405501864	0.0868438700617951\\
10	-15.1579014051165	0.115012379428002	0.124981721117281\\
};
\addplot [color=mycolor4, draw=none, line width=0.6pt, mark=none, mark options={solid, mycolor4}, mark size = 4pt, forget plot]
 plot [error bars/.cd, y dir = both, y explicit]
 table[row sep=crcr, y error plus index=2, y error minus index=3]{%
1	0.443205402663294	0.0157657391566289	0.0159399321312996\\
2	-1.66103997818893	0.0207128481597494	0.0210145605628551\\
3	-3.69201263493529	0.0310795472356711	0.0317638534657791\\
4	-5.51352883566145	0.0476674656814913	0.0492964000556411\\
5	-7.33147619048305	0.0628449140541454	0.0657076640802412\\
6	-9.00985125227136	0.103051617135609	0.110982737950412\\
7	-10.6698356471261	0.0574282865647753	0.0598094023026068\\
8	-12.1958601662628	0.0613471402384658	0.0640720803007842\\
9	-13.6904075895865	0.0819117405501864	0.0868438700617951\\
10	-15.1579014051165	0.115012379428002	0.124981721117281\\
};
\addlegendimage{mycolor4, dashed, line width=0.6pt, mark=*,mark size=1.25pt}
\addlegendentry{$X_\ell^{\bit,\hh}$}

\addplot [color=mycolor3, dashdotted, line width=0.6pt, mark=*, mark size =0.5pt, mark options={solid, mycolor3}, forget plot]
  table[row sep=crcr]{%
1	-0.340038855731138	0.022961359179936	0.0233327214391275\\
2	-2.37632351602757	0.0198892331747813	0.0201672669274329\\
3	-4.36549208822589	0.022283632845105	0.0226332290778792\\
4	-6.23684571271077	0.0230925138233848	0.0234681654185449\\
5	-7.96645353651504	0.0245094510008972	0.0249330402295449\\
6	-9.63589448973896	0.0281166615568882	0.0286755366929459\\
7	-11.2239062692693	0.0302755345990313	0.0309245224442432\\
8	-12.7764948404731	0.0362061120765773	0.0371381883381012\\
9	-14.2557988248202	0.0408585361261853	0.0420495010892772\\
10	-15.7691124058597	0.0432265372428482	0.0445618170864641\\
};
\addplot [color=mycolor3, draw=none, line width=0.6pt, mark=none, mark options={solid, mycolor3}, mark size = 4pt, forget plot]
 plot [error bars/.cd, y dir = both, y explicit]
 table[row sep=crcr, y error plus index=2, y error minus index=3]{%
1	-0.340038855731138	0.022961359179936	0.0233327214391275\\
2	-2.37632351602757	0.0198892331747813	0.0201672669274329\\
3	-4.36549208822589	0.022283632845105	0.0226332290778792\\
4	-6.23684571271077	0.0230925138233848	0.0234681654185449\\
5	-7.96645353651504	0.0245094510008972	0.0249330402295449\\
6	-9.63589448973896	0.0281166615568882	0.0286755366929459\\
7	-11.2239062692693	0.0302755345990313	0.0309245224442432\\
8	-12.7764948404731	0.0362061120765773	0.0371381883381012\\
9	-14.2557988248202	0.0408585361261853	0.0420495010892772\\
10	-15.7691124058597	0.0432265372428482	0.0445618170864641\\
};
\addlegendimage{mycolor3, dashdotted, line width=0.6pt, mark=*,mark size=1.25pt}
\addlegendentry{$X_\ell^{\num,\hh}$}

\end{axis}
\end{tikzpicture}%
\end{subfigure}
}
\caption{Cox-Ingersoll-Ross process at final time: 
bias and variance vs.\ level}
\label{f7}
\end{figure}
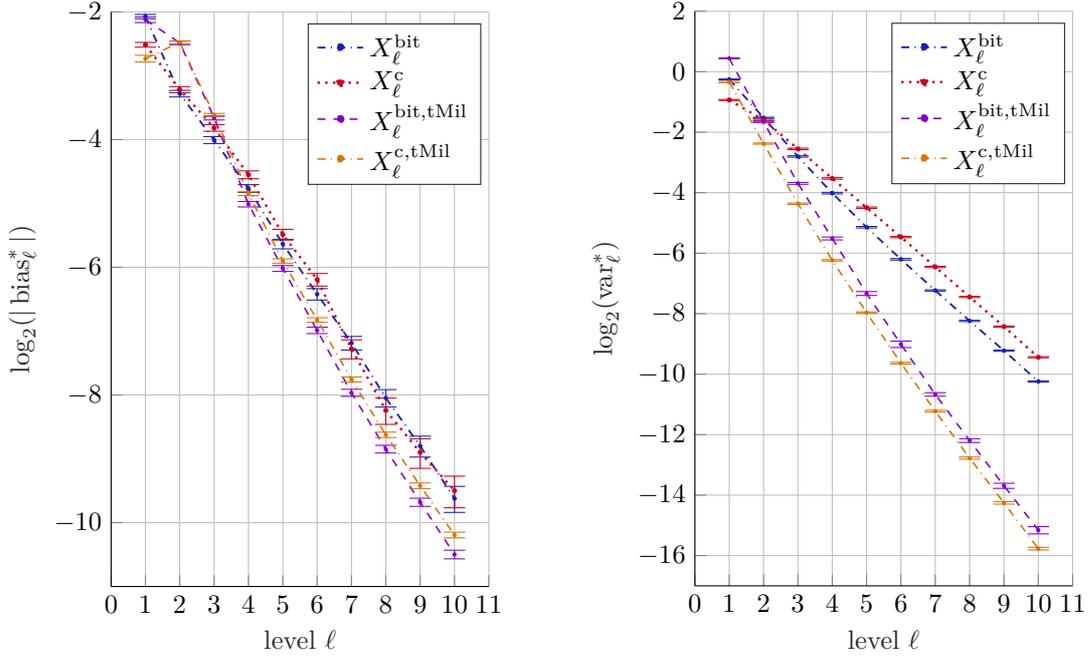

\begin{figure}
\centering
\input{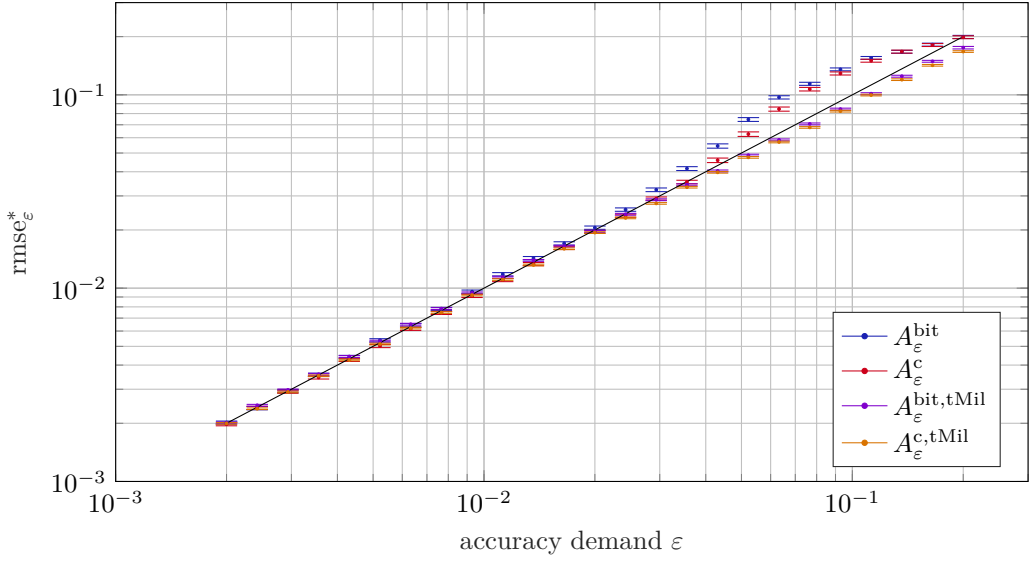}
\caption{Cox-Ingersoll-Ross process at final time: 
root mean squared error vs.\ accuracy demand} 
\label{f8}
\end{figure}

For all four algorithms the root mean squared error is almost equal to 
the accuracy demand, see Figure~\ref{f8},
as is the case of the SDEs considered before.

\begin{figure}
\centering
\input{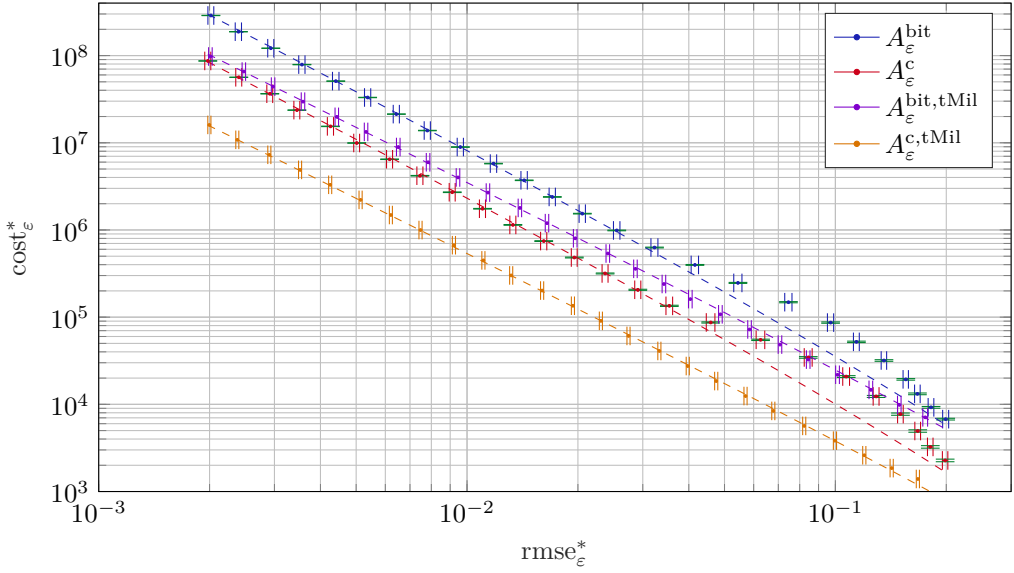}
\caption{Cox-Ingersoll-Ross process at final time: 
cost vs.\ root mean squared error}
\label{f9}
\end{figure}

Finally we relate $\cost^\ast_\eps$ to the root mean squared error 
$\err_\eps^\ast$, see Figure~\ref{f9}.
The exponent $\gamma$ of the logarithmic term is equal to $1.2$ for 
both variants that are based on the Euler scheme
and equal to $0.5$ for both variants that are based on
the truncated Milstein scheme.
The better log-exponent corresponds to the faster decay of the variances.
Furthermore, we have $\kappa^\bit/\kappa^\num=3.51$
and $\kappa^{\bit,\hh}/\kappa^{\num,\hh}=6.52$,
and the multilevel algorithm based on the Euler scheme with random 
numbers has roughly the same root mean squared error as the multilevel 
algorithm based on the truncated Milstein scheme using random bits
in the range considered in Figure~\ref{f9}.

\section*{Acknowledgement}
Mike Giles was partially supported by the UK Engineering and Physical
Science Research Council (EPSRC) through the ICONIC Programme Grant,
EP/P020720/1.
Lukas Mayer was supported by the Deutsche Forschungsgemeinschaft
(DFG) within the RTG 1932 
`Stochastic Models for Innovations in the Engineering Sciences'.

\section*{Appendix}

We present the proofs of Lemma \ref{l1} and Lemma \ref{l2}.
For notational convenience we consider the case $d=1$.
Recall that $Z^{(0,1)},\dots,Z^{(\ell,2^{\ell-1})}$ are independent
and standard normally distributed. 

The distribution of $U^{(i,j)}_\ell$ is symmetric with respect to
$1/2$, so that $\E(Y^{(i,j)}_\ell) = 0$. Moreover,
we have
$\var(Y^{(i,j)}_\ell) = 1$ by construction.
It follows that $\E(V_{k,\ell})=0$ and
\begin{align*}
\var\bigl(V_{k,\ell}\bigr)
&=\bigl(\Delta^{(0,1)}_{k,\ell}\bigr)^2
+
\sum_{i=1}^\ell \sum_{j=1}^{2^{i-1}} 
\bigl(\Delta^{(i,j)}_{k,\ell}\bigr)^2 \\
&=\var\Bigl(\Delta^{(0,1)}_{k,\ell} \cdot Z^{(0,1)} +
\sum_{i=1}^\ell \sum_{j=1}^{2^{i-1}} 
\Delta^{(i,j)}_{k,\ell} \cdot Z^{(i,j)}\Bigr)\\
&=\var\bigl(W_\ell(t_{k,\ell})-W_\ell(t_{k-1,\ell})\bigr)
=2^{-\ell}
\end{align*}
due to \eqref{g15} and the convergence of the L\'evy-Ciesielski
representation.

We have
\[
\var\bigl(\Phi^{-1} \circ U^{(i,j)}_\ell\bigr) \leq 1,
\]
see the end of the proof of \cite[Thm.~1]{GHMR17}, and in
particular for $i=\ell$
\[
\var\bigl(\Phi^{-1} \circ U^{(\ell,j)}_\ell\bigr) \leq 4/5,
\]
which follows from a simple computation. It follows that
\begin{align*}
&\var\bigl(V_{k,\ell}^\prime\bigr) \\
&\qquad =\bigl(\Delta^{(0,1)}_{k,\ell}\bigr)^2
\cdot \var\bigl( \Phi^{-1} \circ U^{(0,1)}_\ell\bigr) 
+
\sum_{i=1}^{\ell} \sum_{j=1}^{2^{i-1}} 
\bigl(\Delta^{(i,j)}_{k,\ell}\bigr)^2 \cdot 
\var\bigl(\Phi^{-1} \circ U^{(i,j)}_\ell\bigr)\\
&\qquad \leq \bigl(\Delta^{(0,1)}_{k,\ell}\bigr)^2
+\sum_{i=1}^{\ell-1} \sum_{j=1}^{2^{i-1}} 
\bigl(\Delta^{(i,j)}_{k,\ell}\bigr)^2 
\mbox{}+\sum_{j=1}^{2^{\ell-1}} 
\bigl(\Delta^{(\ell,j)}_{k,\ell}\bigr)^2 \cdot 
\var\bigl(\Phi^{-1} \circ U^{(\ell,j)}_\ell\bigr)\\
&\qquad \leq 2^{-\ell}
-\sum_{j=1}^{2^{\ell-1}} 
\bigl(\Delta^{(\ell,j)}_{k,\ell}\bigr)^2 \cdot 
\Bigl( 1 -
\var\bigl(\Phi^{-1} \circ U^{(\ell,j)}_\ell\bigr)\Bigr)\\
&\qquad\leq 2^{-\ell}
-1/5\cdot \sum_{j=1}^{2^{\ell-1}} 
\big(\Delta^{(\ell,j)}_{k,\ell}\big)^2
= 9/10 \cdot 2^{-\ell}.
\end{align*}

\bibliographystyle{plain}
\bibliography{bib}

\end{document}